\documentclass[a4paper,11pt]{article}
\usepackage{amsfonts,amssymb,amsmath, dsfont}
\usepackage[english]{babel}
 \makeatletter
 \@addtoreset{equation}{section}
 \makeatother

 \makeatletter
 \@addtoreset{enunciato}{section}
 \makeatother

 \newcounter{enunciato}[section]

 \newtheorem{ittheorem}{Theorem}
 \newtheorem{itlemma}{Lemma}
 \newtheorem{itproposition}{Proposition}
 \newtheorem{itcorollary}{Corollary}
 \newtheorem{itdefinition}{Definition}
 \newtheorem{itremark}{Remark}
 \newtheorem{itclaim}{Claim}
 \newtheorem{itfact}{Fact}
 \newtheorem{itconjecture}{Conjecture}

 \newenvironment{theorem}{\addtocounter{enunciato}{1}
 \begin{ittheorem}}{\end{ittheorem}}

 \newenvironment{lemma}{\addtocounter{enunciato}{1}
 \begin{itlemma}}{\end{itlemma}}

 \newenvironment{proposition}{\addtocounter{enunciato}{1}
 \begin{itproposition}}{\end{itproposition}}

 \newenvironment{corollary}{\addtocounter{enunciato}{1}
 \begin{itcorollary}}{\end{itcorollary}}

 \newenvironment{definition}{\addtocounter{enunciato}{1}
 \begin{itdefinition}}{\end{itdefinition}}

 \newenvironment{remark}{\addtocounter{enunciato}{1}
 \begin{itremark}}{\end{itremark}}

 \newenvironment{claim}{\addtocounter{enunciato}{1}
 \begin{itclaim}}{\end{itclaim}}

 \newenvironment{fact}{\addtocounter{enunciato}{1}
 \begin{itfact}}{\end{itfact}}

 \newenvironment{conjecture}{\addtocounter{enunciato}{1}
 \begin{itconjecture}}{\end{itconjecture}}

 \newcommand{\be}[1]{\begin{equation}\label{#1}}
 \newcommand{\ee}{\end{equation}}

 \newcommand{\bl}[1]{\begin{lemma}\label{#1}}
 \newcommand{\el}{\end{lemma}}

 \newcommand{\br}[1]{\begin{remark}\label{#1}}
 \newcommand{\er}{\end{remark}}

 \newcommand{\bt}[1]{\begin{theorem}\label{#1}}
 \newcommand{\et}{\end{theorem}}

 \newcommand{\bd}[1]{\begin{definition}\label{#1}}
 \newcommand{\ed}{\end{definition}}

 \newcommand{\bcl}[1]{\begin{claim}\label{#1}}
 \newcommand{\ecl}{\end{claim}}

 \newcommand{\bfact}[1]{\begin{fact}\label{#1}}
 \newcommand{\efact}{\end{fact}}

 \newcommand{\bp}[1]{\begin{proposition}\label{#1}}
 \newcommand{\ep}{\end{proposition}}

 \newcommand{\bc}[1]{\begin{corollary}\label{#1}}
 \newcommand{\ec}{\end{corollary}}

 \newcommand{\bcj}[1]{\begin{conjecture}\label{#1}}
 \newcommand{\ecj}{\end{conjecture}}

 \newcommand{\bpr}{\begin{proof}}
 \newcommand{\epr}{\end{proof}}

 \newcommand{\bprlem}[1]{\begin{proofof}{\it Lemma \ref{#1}}.\,\,}
 \newcommand{\eprlem}{\end{proofof}}

 \newcommand{\bprthm}[1]{\begin{proofof}{\it Theorem \ref{#1}}.\,\,}
 \newcommand{\eprthm}{\end{proofof}}

 \newcommand{\bi}{\begin{itemize}}
 \newcommand{\ei}{\end{itemize}}

 \newcommand{\ben}{\begin{enumerate}}
 \newcommand{\een}{\end{enumerate}}

 \newenvironment{proof}{\noindent {\em Proof}.\,\,}{\hspace*{\fill}$\halmos$\medskip}
 \newenvironment{proofof}{\noindent {\em Proof of\,\,}}{\hspace*{\fill}$\halmos$\medskip}
 \newcommand{\halmos}{\rule{1ex}{1.4ex}}

 \newcommand{\one}{{\mathchoice {1\mskip-4mu\mathrm l}
         {1\mskip-4mu\mathrm l}
         {1\mskip-4.5mu\mathrm l}
         {1\mskip-5mu\mathrm l}}}

\def \E {{\mathbb E}}

\def \N {{\mathbb N}}
\def \P {{\mathbb P}}

\def \Z {{\mathbb Z}}

\def \ba {\begin{array}}
\def \ea {\end{array}}

\def\one{\rlap{\mbox{\small\rm 1}}\kern.15em 1}

\begin{document}
\title{Multitype Contact Process on $\mathbb{Z}$: Extinction and Interface}

\author{\renewcommand{\thefootnote}{\arabic{footnote}}
Daniel Valesin
\footnotemark[1]
}

\footnotetext[1]{
Institut de Math\'ematiques, Station 8,
\'Ecole Polytechnique F\'ed\'erale de Lausanne,
CH-1015 Lausanne, Switzerland,
daniel.valesin@epfl.ch
}
\date{March 23, 2010}
\maketitle

\begin{abstract} 
We consider a two-type contact process on $\mathbb{Z}$ in which both types have equal finite range and supercritical infection rate. We show that a given type becomes extinct with probability $1$ if and only if, in the initial configuration, it is confined to a finite interval $[-L, L]$ and the other type occupies infinitely many sites both in $(-\infty, L)$ and $(L, \infty)$. We also show that, starting from the configuration in which all sites in $(-\infty, 0]$ are occupied by type $1$ particles and all sites in $(0, \infty)$ are occupied by type $2$ particles, the process $\rho_t$ defined by the size of the interface area between the two types at time $t$ is tight.
\end{abstract}

\newpage

%%%%%%%%%%%%%% Introduction %%%%%%%%%%%%%%%%%%%%%%%%%%%%%%%%%%%%%%%%%%%%%%%%%%%%%%%%%%%%%%%%%%%%%%%%%%%%%%%%%%%%%%%%%%%%%%

\section{Introduction}
\label{Int}

The contact process on $\Z$ is the spin system with generator
$$\Omega f(\zeta) = \sum_x (f(\zeta^x) - f(\zeta))\; c(x, \zeta);\quad \zeta \in \{0, 1\}^\Z $$
where
$$\left\{\begin{array}{l} \zeta^x(y) = \zeta(y) \text{ if } x \neq y; \\ \zeta^x(x) = 1 - \zeta(x);  \end{array}\right. \quad\qquad c(x, \zeta) = \left\{\begin{array}{cl} 1 &\text{if } \zeta(x) = 1;\\ \lambda \sum_y \zeta(y)\cdot p(y-x)&\text{if } \zeta(x) = 0; \end{array}\right.$$ 
for $\lambda > 0$ and $p(\cdot)$ a probability kernel. We take $p$ to be symmetric and to have finite range $R = \max\{x: p(x) > 0\}$.

The contact process is usually taken as a model for the spread of an infection; configuration $\zeta \in \{0, 1\}^\Z$ is the state in which an infection is present at $x \in \Z$ if and only if $\zeta(x) = 1$. With this in mind, the dynamics may be interpreted as follows: each infected site waits an exponential time of parameter 1, after which it heals, and additionally each infected site waits an exponential time of parameter $\lambda$, after which it chooses, according to the kernel $p$, some other site to which the infection is transmitted if not already present.

We refer the reader to \cite{lig99} for a complete account of the contact process. Here we mention only the most fundamental fact. Let $\bar \zeta$ and $\textbf{0}$ be the configurations identically equal to $1$ and $0$, respectively, $\P_\lambda$ the probability measure under which the process has rate $\lambda$ and $\zeta^0_t$ the configuration at time $t$, started from the configuration where only the origin is infected. There exists $\lambda_c$, depending on $p$, such that\\
$\bullet\;$ if $\lambda \leq \lambda_c$, then $\P_\lambda(\zeta^0_t \neq \textbf{0} \; \forall t) = 0$ and $\delta_{\bar \zeta}S(t) \to \delta_{\textbf{0}}$;\\
$\bullet\;$ if $\lambda > \lambda_c$, then $\P_\lambda(\zeta^0_t \neq \textbf{0} \; \forall t) > 0$ and $\delta_{\bar \zeta}S(t)$ converges, as $t \to \infty$, to some non-trivial invariant measure.\\
Again, see \cite{lig99} for the proof. Throughout this paper, we fix $\lambda > \lambda_c$.

The multitype contact process was introduced in \cite{neuhauser} as a modification of the above system. Here we consider a two-type contact process, defined as the particle system $(\xi_t)_{t \geq 0}$ with state space $\{0, 1, 2\}^\Z$ and generator
$$\begin{aligned}
\Lambda f(\xi) = &\sum_{x: \xi(x) \neq 0} (f(\xi^{x, 0})-f(\xi)) + \\&\sum_{x: \xi(x) = 0} \big[(f(\xi^{x,1})- f(\xi)) \; c_1(x, \xi) + (f(\xi^{x,2})- f(\xi)) \; c_2(x, \xi)\big]; \quad \xi \in \{0, 1, 2\}^\Z,
\end{aligned}$$
where
$$\begin{array}{l}\left\{ \begin{array}{l} \xi^{x, i}(y) = \xi(y) \text{ if } x \neq y; \\ \xi^{x, i}(x) = i, \end{array} \right.\\\quad i=0,1,2;\end{array}
\qquad \begin{array}{l} c_i(x, \xi) =  \lambda \sum_y \mathds{1}_{\{\xi(y) = i\}} \cdot p(y - x),\medskip \\\quad i = 1, 2. \end{array}
$$
($\mathds{1}$ denotes indicator function).

This is thought of as a model for competition of two biological species. Each site in $\Z$ corresponds to a region of space, which can be either empty or occupied by an individual of one of the two species. Occupied regions get empty at rate $1$, meaning natural death of the occupant, and empty regions get occupied by a rate that depends on the number of individuals of each species living in neighboring sites, and this means a new birth. The important point is that occupancy is strong in the sense that, if a site has an individual of, say, type $1$, the only way it will later contain an individual of type $2$ is if the current individual dies and a new birth occurs originated from a type $2$ individual.

Let us point out some properties of the above dynamics. First, it is symmetric for the two species: both die and give birth following the same rules and restrictions. Second, if only one of the two species is present in the initial configuration, then the process evolves exactly like in the one-type contact process. Third, if we only distinguish occupied sites from non-occupied ones, thus ignoring which of the two types is present at each site, again we see the evolution of the one-type contact process.

The first question we address is: for which initial configurations a given type (say, type $1$) becomes extinct with probability one? By extinction we mean: for some time $t_0$ (and hence all $t \geq t_0$), $\xi_{t_0}(x) \neq 1$ for all $x$. We prove
\bt{extThm}
Assume at least one site is occupied by a $1$ in $\xi_0$. The $1$'s become extinct with probability one if and only if there exists $L > 0$ such that\\
$i.)\;\xi_0(x) \neq 1 \; \forall x \notin [-L, L]$ and\\
$ii.) \; \#\{x \in (-\infty, -L]: \xi_0(x) = 2\} = \# \{x \in [L, \infty): \xi_0(x) = 2\} = \infty$\et
($\#$ denotes cardinality). This result is a generalization of Theorem 1.1. in \cite{amp}, which is the exact same statement in the nearest neighbour context (\textit{i.e.}, $p(1) = p(-1) = 1/2$). Althought there are some points in common between our proof and the one in that work, our general approach is completely different.

Now assume that the range $R > 1$. Define the ``heaviside'' configuration as $\xi^H = \mathds{1}_{(-\infty, 0]} + 2 \cdot \mathds{1}_{(0, \infty)}$ and denote by $\xi_t$ the two-type contact process with initial condition $\xi_0 = \xi^H$. Define 
$$r_t = \sup\{x: \xi_t(x) = 1\}, \quad l_t = \inf\{x: \xi_t(x) = 2\}, \quad \rho_t = r_t - l_t.$$ 
We have $\rho_0 = -1$, and at a given time $t$ both events $\{\rho_t > 0\}$ and $\{\rho_t < 0\}$ have positive probability. If $\rho_t > 0$, we call the interval $[l_t, r_t]$ the \textit{interface area}. The question we want to ask is: if $t$ is large, is it reasonable to expect a large interface? We answer this question negatively.
\bt{interThm}
The law of $(\rho_t)_{t \geq 0}$ is tight; that is, for any $\epsilon > 0$, there exists $L > 0$ such that $\P(|\rho_t| > L) < \epsilon$ for every $t \geq 0$.
\et

There are several works concerning interface tightness in one-dimensional particle systems, the first of which is \cite{cd}, where interface tightness is established for the voter model. Others are \cite{samir1}, \cite{samir2}, \cite{swart} and  \cite{ampv}. 

In \cite{ampv}, it is shown that interface tightness also occurs on another variant of the contact process, namely the grass-bushes-trees model considered in \cite{durswin}, with both species having same infection rate and non-nearest neighbor interaction. The difference between the grass-bushes- trees model and the multitype contact process considered here is that, in the former, one of the two species, say the 1's, is privileged in the sense that it is allowed to invade sites occupied by the 2's. For this reason, from the point of view of the 1's, the presence of the 2's is irrelevant. It is thus possible to restrict attention to the evolution of the 1's, and it is shown that they form barriers that prevent entrance from outside; with this at hand, interface tightness is guaranteed irregardless of the evolution of the 2's. Here, however, we do not have this advantage, since we cannot study the evolution of any of the species while ignoring the other.

Our results depend on a careful examination of the temporal dual process; that is, rather than moving forward in time and following the descendancy of individuals, we move backwards in time and trace ancestries. The dual of the multitype contact process was first studied by Neuhauser in \cite{neuhauser} and may be briefly described as follows. Each site $x \in \Z$ at (primal) time $s$ has a random (and possibly empty) \textit{ancestor sequence}, which is a list of sites $y \in \Z$ such that an infection could potentially be transmitted from $(y, 0)$ to $(x, s)$. The ancestors on the list are ranked in decreasing order; the idea is that if the first ancestor is not occupied in $\xi_0$, then we look at the second, and so on, until we find the first on the list that is occupied in $\xi_0$, and take its type as the one passed to $x$. We denote this sequence $(\eta^x_{1, s}, \eta^x_{2, s}, \ldots)$. By moving in time in the opposite sense as that of the original process and using the graphical representation of the contact process for ``negative'' primal times, we can define the \textit{ancestry process of} $x$, $((\eta^x_{1, t}, \eta^x_{2, t}, \ldots))_{t \geq 0}$. The process given by the first element of the sequence, $(\eta^x_{1, t})_{t \geq 0}$, is called the \textit{first ancestor} process. We point out three key properties of the ancestry process:

\begin{itemize}
\item \textbf{First ancestors behave as random walks.} In \cite{neuhauser} it is proven that, on the event that a site $x$ has a nonempty ancestry at all times $t \geq 0$, we can define an increasing sequence of random \textit{renewal times} $(\tau^x_n)_{n \geq 0}$ with the property that the space-time increments $(\eta^x_{1, \tau^x_{n+1}}-\eta^x_{1, \tau^x_n}, \tau^x_{n+1} - \tau^x_n)$ are independent and identically distributed. This fact enormously simplifies the study of the first ancestor process, which is not markovian and at first seems very complicated.

\item \textbf{Ancestries coalesce.} If we are to use the dual process to obtain information about the joint distribution of the states at sites $x$ and $y$ at a given time, we must study the joint behavior of two ancestry processes, specially of two first ancestor processes. The intuitive picture is that this behavior resembles that of two random walks that are independent until they meet, at which time they coalesce. We give a new approach to formalize this notion, one that we believe provides a clear understanding of the picture and allows for detailed results.

In order to follow two first ancestor processes simultaneously, we define \textit{joint renewals} $(\tau^{x,y}_n)_{n \geq 0}$ and argue that the law of the processes after a joint renewal only depends on their initial difference at the instant of the renewal. Thus, the discrete-time process defined by the difference between the two processes at the instants of renewals is a Markov chain on $\Z$. For this chain, zero is an absorbing state and corresponds to coalescence of first ancestors. We also show that, far from the origin, the transition probabilities of the chain become close to a symmetric measure on $\Z$, and from this fact we are able to show that the tail of the distribution of the hitting time of $0$ for the chain looks like the one associated to a simple random walk on $\Z$. From this construction and estimate we also bound the expected distance between ancestors at a given time.
\item \textbf{Ancestries become sparse with time.} Consider the system of coalescing random walks in which each site of $\Z$ starts with one particle at time $0$. The density of occupied sites at time $t$, which is equal to the probability of the origin being occupied, tends to $0$ as $t \to \infty$. We prove a similar result for our ancestry sequences. Fix a \textit{truncation level} $N$ and, at dual time $t$, mark the $N$ first ancestors of each site at dual time $0$ (this gives the set $\{\eta^x_{n, t}: 1 \leq n \leq N, x \in \mathbb{Z}: \text{ the ancestry of } x \text{ reaches time } t\}$). We show that the density of this random set tends to $0$ as $t \to \infty$, and estimate the speed of this convergence depending on $N$.
\end{itemize}

From this last fact, we can immediately prove Theorem \ref{extThm} under the stronger hypothesis that all sites outside $[-L, L]$ are occupied by $2$'s in $\xi_0$. To obtain the general case, we then use a structure called a \textit{descendancy barrier}, whose existence was established in \cite{ampv}.

The proof of Theorem \ref{interThm} is more intricate. It follows the main steps of the argument in \cite{cd} for the voter model. Starting from $\xi^H$, say that sites $x < y$ form a $k$\textit{-inversion} at (primal) time $t$ if $y - x = k,\; \xi_t(x) = 2$ and $\xi_t(y) = 1$. Using the coalescence properties described above, it is shown that the expected number of $k$-inversions at time $t$ is uniformly bounded in $k$ and $t$. A consequence is that, if instead of looking at the whole configuration $\xi_t$, we only look at its restriction to a sparse subset $R \subset \Z$, it is unlikely (uniformly in $t$) that we find any inversion at all. Next, given $0 < s < t$, consider the random set $R(s, t)$ of sites that are occupied in $\xi_s$ and that survive up to time $t$ (that is, $x \in R(s,t)$ if there exists an ``infection path'' starting at some site at time $0$, passing through $(x, s)$ and reaching time $t$ such that every jump in this path lands on an unoccupied site). We show that, if $t - s$ is large, then the density of $R(s, t)$ is small, and this is uniform in $t$. Putting these facts together, we conclude that, for $0 < s < t$ appropriately chosen, with large probability no inversion is present in the restriction of $\xi_s$ to $R(s, t)$. This means that the inversions that are present in $\xi_t$ are formed in the final time interval $[t - s, t]$. Fixing large $s$ and changing $t$, we get tightness of the number of inversions, and it is then straightforward to establish tightness of the interface size.

We believe that our results and general approach may prove useful in other questions concerning the multitype contact process, in particular those that relate to properties of trajectories, of which not much is known.

The author would like to thank Thomas Mountford for all his help and the colleagues and friends Augusto Teixeira, Johel Beltr\'an and Renato Santos for helpful ideas and encouragement.

%%%%%%%%%%%%%% section 1 %%%%%%%%%%%%%%%%%%%%%%%%%%%%%%%%%%%%%%%%%%%%%%%%%%%%%%%%%%%%%%%%%%%%%%%%%%%%%%%%%%%%%%%%%%%%%%

\section{Ancestry process}
\label{Anc}

We will start describing the familiar construction of the one-type contact process from its graphical representation. We will then show how the same representation can be used to construct the multitype contact process, present the definition of the ancestry process together with some facts from \cite{neuhauser}, and finally prove a  simple lemma.

Suppose given a collection of independent Poisson processes on $[0, \infty)$: 
$$(D^x)_{x \in \Z} \text{ of rate } 1, \quad (N^{(x, y)})_{x, y \in \Z} \text{ of rate } \lambda \cdot p(y - x).$$
A \textit{Harris construction} $H$ is a realization of all such processes. $H$ can thus be understood as a point measure on $(\Z \cup \Z^2) \times [0, \infty)$. Sometimes we abuse notation and denote the collection of processes itself by $H$. Given $(x, t) \in \Z \times [0, \infty)$, let $\theta(x, t)(H)$ be the Harris construction obtained by shifting $H$ so that $(x, t)$ becomes the space-time origin. By translation invariance of the space-time construction, $\theta(x, t)(H)$ and $H$ have the same distribution. We will also write $H_{[0, t]}$ to denote the restriction of $H$ to $\Z \times [0, t]$, and refer to such restrictions as finite-time Harris constructions.

Given a Harris construction $H$ and $(x, s), (y, t) \in \Z \times [0, \infty)$ with $s < t$, we write $(x, s) \leftrightarrow (y, t)$ (in $H$) if there exists a piecewise constant, right-continuous function $\gamma: [s, t] \to \Z$ such that\\
$\bullet\;\gamma(s) = x, \gamma(t) = y;$\\
$\bullet\;\gamma(r) \neq \gamma(r-)$ if and only if $r \in N^{(\gamma(r-), \gamma(r))}$;\\
$\bullet\; \nexists s \leq r \leq t$ with $r \in D^{\gamma(r)}$.\\
One such function $\gamma$ is called a \textit{path} determined by $H$. The points in the processes $\{D^x\}$ are usually called \textit{death marks}, and the points in $\{N^{(x,y)}\}$ are called \textit{arrows}. Thus, a path can be thought of as a line going up from $(x, s)$ to $(y, t)$ following the arrows and not crossing any death marks.

Given $A \subset \Z, (x, t) \in \Z \times [0, \infty)$ and a Harris construction $H$, put $$[\zeta^A_t(x)](H) = \mathds{1}_{\{\text{For some } y \in A, (y, 0) \leftrightarrow (x, t) \text{ in } H\}}.$$ Under the law of $H$, $(\zeta^A_t)$ has the distribution of the contact process with parameter $\lambda$, kernel $p$ and initial state $\mathds{1}_A$; see \cite{d1} for details. From now on, we omit dependency on the Harris construction and write (for instance) $\zeta_t$ instead of $\zeta_t(H)$.

Before going into the multitype contact process, we list some properties of the one-type contact process that will be very useful. Fix $(x,s) \in \Z\times [0, \infty)$ and $t > s$. Define the time of death and maximal distance traveled until time $t$ for an infection that starts at $(x, s)$,
$$\hat T^{(x, s)} = \inf\{s' > s: \nexists y: (x, s) \leftrightarrow (y, s'),$$
$$M^{(x, s)}_t = \sup\{|y - x|: (x, s) \leftrightarrow (y, s') \text{ for some } s' \in [s, t]\}$$
(these only depend on $H$ and are thus well-defined irregardless of $\xi_s(x)$).
When $s = 0$, we omit it and write $\hat T^x, M^x_t$. If $A \subset \Z$, we also define $\hat T^A = \inf\{t \geq 0: \nexists\;  x \in A, y \in \Z: (x, 0) \leftrightarrow (y, t)\}$. We start mentioning that $M^{(x, s)}_t$ is stochastically dominated by a sum of Poisson random variables, so there exist $\kappa , c, C > 0$ such that
\be{comparePoi}\P(M^x_t > \kappa t) \leq Ce^{-ct} \quad \forall x \in \Z, t \geq 0. \ee
Next, since we are taking $\lambda > \lambda_c$, we have $\P(\hat T^x = \infty) = \P(\hat T^0 = \infty) > 0$ for all $x$, and
\be{poscorr}\P(\hat T^x = \hat T^y = \infty) \geq \P(\hat T^0 = \infty)^2 > 0, \qquad \forall x, y \in \Z.\ee
This follows from the self-duality of the contact process and the fact that its upper invariant measure has positive correlations; see \cite{lig85}. Our last property is that there exist $c, C > 0$ such that, for any $A \subset \Z$ and $t > 0$,
\be{fastDeath} \P(t < \hat T^A < \infty) \leq Ce^{-ct}. \ee
For the case $R = 1$, this is Theorem 2.30 in \cite{lig99}. The proof uses a comparison with oriented percolation and can be easily adapted to the case $R > 1$.

To obtain the graphical representation for the multitype contact process, we have to proceed as above but ignore the arrows whose endpoints are already occupied. This was first done in \cite{neuhauser}; there, an algorithmic procedure is provided to find the state of each site at a given time. Here we provide an approach that is formally different but amounts to the same. Fix $(x, t) \in \Z \times [0, \infty)$, a Harris construction $H$ and $\xi_0 \in \{0, 1, 2\}^\Z$. Let $\Gamma$ be the set of paths $\gamma$ that connect points of $\Z \times \{0\}$ to $(x, t)$ in $H$. Assume that $\#\Gamma < \infty$; this happens with probability one if $H$ is sampled from the processes described above. For the moment, also assume that $\Gamma \neq \emptyset$. Given $\gamma, \gamma ' \in \Gamma$, let us write $\gamma \prec \gamma '$ if there exists $\bar s \in (0, t)$ such that $\gamma(s) = \gamma '(s) \; \forall s \in [\bar s, t]$ and $\gamma(\bar s) \neq \gamma(\bar s-), \gamma '(\bar s) = \gamma '(\bar s-)$. From the fact that these paths are all piecewise constant, have finitely many jumps and the same endpoint, we deduce that $\prec$ is a total order on $\Gamma$. We can then find $\gamma^*_0$, the maximal path in $\Gamma$. Now define $\Gamma_1 = \{\gamma \in \Gamma: \gamma(0) \neq \gamma^*_0(0)\}$ and $\gamma^*_1$ as the maximal path in $\Gamma_1$. Then define $\Gamma_2 = \{\gamma \in \Gamma_1: \gamma(0) \neq \gamma^*_1(0)\}$, and so on, until $\Gamma_N = \emptyset$. For $0 \leq n < N$, denote $\hat \eta^x_{n, t} = \gamma^*_n(0)$, and for $n \geq N$ put $\hat \eta^x_{n, t} = \triangle$. We claim that
\be{obstruct} \forall n < N, \forall s \text{ such that } \gamma^*_n(s-) \neq \gamma^*_n(s), \text{ we have } (\gamma^*_n(s), s) \notin \zeta^{\{\hat \eta^x_{n+1, t}, \ldots, \hat \eta^x_{N-1, t}\}}_s\ee
(Here $\zeta_\cdot$ continues to denote the one-type contact process defined from $H$). In words, each jump of $\gamma^*_n$ lands on a space-time point that cannot be reached by paths coming from $\hat \eta^x_{n+1, t}, \ldots, \hat \eta^x_{N-1, t}$. If this were not the case, we could obtain $m < n, s \in [0, t]$ and $\gamma$ with $\gamma(0) = \hat \eta^x_{n, t}$ and $\gamma^*_m(s-) \neq \gamma^*_m(s) = \gamma(s)$. But we could then construct a path $\gamma '$ coinciding with $\gamma$ on $[0, s]$ and with $\gamma^*_m$ on $(s, t]$, and $\gamma '$ would contradict the maximality that defined $\gamma^*_m$.

If $\xi_0(\hat \eta^x_{n, t}) = 0$ for all $n < N$, put $\xi_t(x) = 0$. Otherwise, if $k= \min\{n: \xi_0(\hat \eta^x_{n, t}) \neq 0\}$, put $\xi_t(x) = \xi_0(\hat \eta^x_{k, t})$. In this second case, using (\ref{obstruct}), we see that there is a path connecting $(\hat \eta^x_{k, t}, 0)$ to $(x, t)$ which obstructs all paths connecting $\{y \neq \hat \eta^x_{k, t}: \xi_0(y) \neq 0\}  \times \{0\}$ to $(x, t)$ and is not obstructed by any of them. Finally, if $\Gamma = \emptyset$, put $\hat \eta^x_{n, t} = \triangle$ for every $n$ and set $\xi_t(x) = 0$. It now follows that $(\xi_t(x))_{x \in \Z}$ has the distribution of the multitype contact process at time $t$, with initial state $\xi_0$.

By considering the time dual of the above construction, we will now define the ancestry process, our main object of investigation. Again fix $x \in \Z, t > 0$ and a Harris construction $H = ((D^x), (N^{(x, y)}))$. Let $\mathcal{I}_t(H)$ be the finite-time Harris construction on $[0, t]$ obtained from $H_{[0, t]}$ by inverting the direction of time and of the arrows; formally, $\mathcal{I}_t(H) = ((\hat D^x), (\hat N^{(x, y)}))$, where
$$\hat D^x(s) = D^x(t-s), \quad \hat N^{(x, y)}(s) = N^{(y,x)}(t-s), \quad 0\leq s \leq t, x, y \in \Z.$$
Two immediate facts are that the laws of $H_{[0, t]}$ and $\mathcal{I}_t(H)$ are equal and that $(x, 0) \leftrightarrow (y, t)$ in $H$ if and only if $(y, 0) \leftrightarrow (x, t)$ in $\mathcal{I}_t(H)$. Define $\eta^x_{n, t}(H) = \hat \eta^x_{n, t}(\mathcal{I}_t(H))$. The $(\Z \cup \{\triangle\})^\infty$-valued process
$$t \mapsto (\eta^x_{1, t}, \eta^x_{2, t}, \ldots)$$
is called the \textit{ancestry process} of $x$. $\eta^x_{n, t}$ is called the $n$\textit{th ancestor of} $x$ \textit{at time} $t$. As a repetition of what was stated in the last paragraph, if we have $\xi_0 \in \{0, 1, 2\}^\Z, t > 0$ and the sequences $(\eta^x_{n,t})_{n \geq 1}$ for each $x \in \Z$, then we can define
$$\xi_t(x) = \left\{\begin{array}{ll}0, &\text{if for each } n, \text{ either } \xi_0(\eta^x_{n, t}) = 0 \text{ or } \eta^x_{n, t} = \triangle;\\ \xi_0(\eta^x_{n^*(x), t}), &\text{where } n^*(x) = \inf\{n: \xi_0(\eta^x_{n, t}) \neq 0\},\end{array}\right.$$
and then $\xi_t$ has the law of the multitype contact process at time $t$ started from $\xi_0$.

We will employ the expressions ``primal time" and ``dual time" referring to the evolution of the original process $t \mapsto \xi_t$ and of the ancestry process $t \mapsto (\eta^x_{n,t})_{n \geq 1}$ respectively; of course, it only makes sense to consider both processes simultaneously if we fix some time $t$, place the primal time origin at $t$ and think of primal time as decreasing from $t$ to $0$ as dual time increases from $0$ to $t$. However, the definition in the previous paragraph allows us to obtain the ancestry process from a Harris construction $H$ for all positive times. From now on, unless explicitly stated otherwise, whenever we mention the Harris construction $H$ and functions of it, such as $\hat T^{(x,s)}$ and $M^{(x, s)}_t$, we mean the Harris construction used to define the ancestry process.

Given $x \in \Z$ and $0 \leq s \leq t$, we define $\eta^{(x,s)}_{n, t} = \eta^x_{n, t-s}(\theta(0, s)(H))$ (that is, the $n$th ancestor in the graph that grows from $(x, s)$ up to time $t$). Also, when $n = 1$, we omit it, writing $\eta^x_t, \eta^{(x, s)}_t$ instead of $\eta^x_{1, t}, \eta^{(x, s)}_{1, t}$. Finally, we write $\eta^{(x,s)}_{*, t} = \{\eta^{(x, s)}_{n, t} \in \Z: n \geq 1\}$, and similarly for $\eta^x_{*, t}$.

The following is an easy consequence of the definition of the ancestry process with the ordering of paths $\prec$ defined above.
\bl{subdesc}
$(i.)$ Let $s > 0$, assume that $\eta^x_s \neq \triangle$ and $\hat T^{(\eta^x_s, s)} = \infty$. Then, for every $t \geq s, \eta^x_t = \eta^{(\eta^x_s, s)}_t$.\\
$(ii.)$ Let $0 \leq s < t, z_1, \ldots, z_N \in \Z$ and assume
\be{desctrans}
\begin{array}{c} 
\eta^x_{i, s} \neq \triangle, \;\eta^{(\eta^x_{i,s}, s)}_{*, t} = \emptyset, \quad 1 \leq i < n \medskip \\
\eta^x_{n,s} \neq \triangle, (\eta^{(\eta^x_{n,s}, s)}_{1, t}, \ldots, \eta^{(\eta^x_{n,s}, s)}_{N, t})=(z_1, \ldots, z_N) \end{array} \nonumber \ee
(that is, the first $n-1$ ancestors of $x$ at time $s$ do not reach time $t$, but the $n$-th one does, with ancestors $z_1, \ldots, z_N$). Then, $$(\eta^x_{1, t}, \ldots, \eta^x_{N, t}) = (z_1, \ldots, z_n).$$
\el

Given $x \in \Z$, on $\{\hat T^x = \infty\}$, define $\tau^x_0 \equiv 0$,
$$\tau^x_1 = \inf\{t \geq 1: \hat T^{(\eta^x_t, t)} = \infty\},$$
and, for $n \geq 1$, on $\{\hat T^x = \infty, \tau^x_n < \infty, \eta^x_{\tau^x_n} = z\}$, define
$$\tau^x_{n+1} = \tau^x_n + \tau^0_1 \circ \theta(z, \tau^x_n).$$
For the sake of readability, we will sometimes write $\tilde \P^x(\cdot)$ and $\tilde \E^x(\cdot)$ instead of $\P(\cdot|\hat T^x = \infty)$ and $\E(\cdot|\hat T^x = \infty)$. 

In \cite{neuhauser}, it is shown that under $\tilde \P^x$, the times $\tau^x_n$ work as renewal times for the process $\eta^x_t$, that is, the (Time length, Trajectory) pairs 
$$(\tau^x_{n+1} - \tau^x_n,\; t \in [0, \tau^x_{n+1} - \tau^x_n] \mapsto \eta^x_{\tau^x_n + t} - \eta^x_{\tau^x_n})$$
are independent and identically distributed. This follows from an idea of Kuczek (\cite{kuk}) which by now is an important tool in the particle systems literature. In our current setting, it can be explained as follows. The probability $\tilde \P^x$ is the original probability for the process conditioned on the event $\{(x, 0) \text{ lives forever}\}$. But $(x, 0)$ being connected to $(\eta^x_{\tau^x_1}, \tau^x_1)$ and $(\eta^x_{\tau^x_1}, \tau^x_1)$ living forever imply that $(x, 0)$ lives forever, the event of the former conditioning. This and the fact that, under $\P$, restrictions of $H$ to disjoint time intervals are independent yield that, under $\tilde \P^x$, the shifted Harris construction $\theta(\eta^x_{\tau^x_1}, \tau^x_1)(H)$ has same law as $H$. The argument is then repeated for all $\tau^x_n, n\geq 1$. 

The Proposition below lists the properties of the renewal times that we will need. The proof is in \cite{neuhauser}, except for part $(ii.)$, which is an adaption of Lemma 7 in \cite{msweet} to our context.

\bp{renewal} $(i.)\; \tilde \P^0(\tau_n^0 < \infty) = 1 \; \forall n.$\medskip\\
$(ii.)\;$ For $n \geq 0$, let
$$H_n = H_{[0, \tau^0_n(H)]}, \qquad H_{n+} = \theta(\eta^0_{\tau^0_n}, \tau^0_n)(H).$$
Given an event $A$ on finite-time Harris constructions and an event $B$ on Harris constructions, we have
$$\tilde \P^0(H_n \in A , H_{n+} \in B) = \tilde \P^0(H_n \in A) \cdot \tilde \P^0(H \in B).$$
$(iii.)\;$ Under $\tilde \P^x$, the $\Z$-valued process $\big(\eta^x_{\tau^x_n}\big)_{n \geq 0}$ is a symmetric random walk starting at $x$ and with transitions
$$P(z, w) = \tilde \P^0\big(\eta^0_{\tau_1^0} = w-z\big).$$
$(iv.) \;$ There exist $c, C > 0$ such that
$$\tilde \P^0\big(\tau_1^0 \vee M_{\tau_1}^0 > r\big) \leq Ce^{-cr}.$$
\ep

To conclude this section, we prove some simples properties of the first ancestor process.

\br{triangle} Every time we write events involving a random variable $\eta$ that may take the value $\triangle$, such as $\{\eta \leq 0\}$, we mean $\{\eta \neq \triangle, \eta \leq 0\}$. This applies to part $(iii)$ of the following lemma. Also, we convention to put $\E(f(\eta)) = \E(f(\eta); \eta \neq \triangle)$ for every function $f$. \er

\bl{exp0}$(i.)$ There exist $c, C > 0$ such that, for all $0 \leq a < b$,
$$\tilde \P^0(\nexists n: \tau_n^0 \in [a, b]) \leq Ce^{-c(b-a)}.$$
$(ii.)$ There exists $C > 0$ such that, for all $0 \leq s < t$,
$$\tilde \E^{0} \big(\; (\eta_t^0)^2 - (\eta_s^0)^2 \;\big) \leq C + C(t-s).$$
$(iii.)$ There exist $c, C > 0$ such that for all $l \geq 0$,
$$\P(|\eta^0_t| > l) \leq Ce^{-cl^2/[t]} + Ce^{-cl}.$$
\el
\bpr Define on $\{\hat T^0 = \infty\}$, for $t \geq 0$,
\be{defpsi} \tau_{t-} = \sup\{\tau^0_n \leq t: n \in \N\},\qquad \tau_{t+} = \inf\{\tau^0_n \geq t: n \in \N\},\qquad \psi_t = M^{\left(\eta_{\tau_{t-}}, \tau_{t-}\right)}_{\tau_{t+}} \vee (\tau_{t+} - \tau_{t-}). \nonumber \ee
Using Proposition \ref{renewal}$(ii.)$ and $(iv.)$,
\begin{eqnarray}
&&\tilde \P^0(\psi_t > x) = \sum_{k=0}^\infty \tilde \P^0(\tau^0_k < t,\; \tau^0_{k+1} \geq t,\; \psi_t > x) \nonumber \\
&& \; = \sum_{k=0}^\infty \int_0^t \tilde \P^0\left(\tau^0_1 \geq t-s, \;  M^0_{\tau^0_1} \vee \tau^0_1 > x \right) \; \tilde \P^0(\tau^0_k \in ds) \nonumber \\
&& \; \leq \sum_{k=0}^\infty \sum_{i=1}^{\lceil t \rceil} \int_{i-1}^i\left[ \tilde \P^0\left(\tau^0_1 \geq t-s\right) \wedge \tilde \P^0\left(M^0_{\tau^0_1} \vee \tau^0_1 > x \right)\right] \; \tilde \P^0(\tau^0_k \in ds) \nonumber \\
&& \; \leq \sum_{i=1}^{\lceil t \rceil} \left[Ce^{-c(t-i)} \wedge Ce^{-cx}\right] \; \sum_{k=0}^\infty \tilde \P^0\left(\tau^0_k \in [i-1, i]\right) \nonumber \\
&& \; = \sum_{i=1}^{\lceil t \rceil} \left[Ce^{-c(t-i)} \wedge Ce^{-cx}\right] \; \tilde \E^0|\{n: \tau^0_n \in [i-1, i]\}| . \label{withexp}
\end{eqnarray}
Observe that the above expectation is less than $1$, because there is at most one renewal in each unit interval. (\ref{withexp}) is thus less than $$C\sum_{i=1}^\infty [e^{-ci} \wedge e^{-cx}] \leq C \lceil x \rceil e^{-cx} + C \sum_{i={\lceil x \rceil + 1}}^\infty e^{-ci} \leq C e^{-cx};$$
since this does not depend on $t$, we get
\be{boundpsitail} \tilde \P(\psi_t > x) \leq Ce^{-cx} \ee
for some $c, C > 0$ and all $t \geq 0$. Let us now prove the two statements of the Lemma.

$(i.)$ For $0 \leq a < b$,
$$\tilde \P^0(\nexists n: \tau^0_n \in [a, b]) = \tilde \P^0(\tau_{a+}-\tau_{a-} > b-a) \leq \tilde \P^0(\psi_a > b-a) \leq Ce^{-c(b-a)}.$$

$(ii.)$ The definition of $\psi_t$ and (\ref{boundpsitail}) imply
\begin{eqnarray}
&&|\eta^0_t - \eta^0_{\tau_{t-}}|, |\eta^0_{\tau_{t+}} - \eta^0_t|, |\eta^0_{\tau_{t+}} - \eta^0_{\tau_{t-}}| \leq 2 \psi_t; \label{comparepsi}\\
&&\qquad \quad\sup_{t \geq 0} \;\tilde \E^0((\psi_t)^2) < \infty. \label{boundpsi}\end{eqnarray}

Next, note that, for any $t > 0$,
$$\tilde \E^0\left((\eta^0_{\tau_{t+}})^2\right)  \leq \tilde \E^0\left(\max_{1 \leq i \leq \lceil t \rceil}\;(\eta^0_{\tau^0_i})^2\right)$$
since there are at most $\lceil t \rceil$ renewals in $[0, t]$. 
By the reflection principle (see \cite{d2}, page 285), the expectation on the right-hand side is less than $2\; \tilde \E^0((\eta^0_{\tau^0_{\lceil t \rceil}})^2) = 2 \;t\; \tilde \E^0((\eta^0_{\tau^0_1})^2)$, so we have
\be{boundtaut}\tilde \E^0((\eta^0_{\tau_{t+}})^2) \leq C \cdot t.\ee

With (\ref{comparepsi}), (\ref{boundpsi}) and (\ref{boundtaut}) at hand, we are ready to estimate
\begin{eqnarray}&&\tilde \E^0\left(\;(\eta^0_t)^2 - (\eta^0_s)^2 \;\right) =\nonumber\\
&&\qquad\tilde \E^0\left(\;(\eta^0_t)^2 - (\eta^0_{\tau_{t+}})^2 \;\right)+ \tilde \E^0\left(\;(\eta^0_{\tau_{s+}})^2 - (\eta^0_s)^2 \;\right) + \tilde \E^0\left(\;(\eta^0_{\tau_{t+}})^2 - (\eta^0_{\tau_{s+}})^2 \;\right). \label{threeterms}
\end{eqnarray}
Let us treat each of the three terms separately. Using the independence of increments between different pairs of renewals and (\ref{comparepsi}), we have
\begin{eqnarray} 
&& \tilde \E^0\left((\eta^0_t)^2 - (\eta^0_{\tau_{t+}})^2\right) = \nonumber\\
&&\tilde \E^0\big((\eta^0_{\tau_{t-}})^2 + (\eta^0_t - \eta^0_{\tau_{t-}})^2 + 2\eta^0_{\tau_{t-}}(\eta^0_t - \eta^0_{\tau_{t-}}) -\nonumber\\
&&\qquad\qquad\qquad\qquad (\eta^0_{\tau_{t-}})^2 - (\eta^0_{\tau_{t+}} - \eta^0_{\tau_{t-}})^2 - 2\eta^0_{\tau_{t-}}(\eta^0_{\tau_{t+}} - \eta^0_{\tau_{t-}})\big) = \nonumber \\
&&\tilde \E^0((\eta^0_t - \eta^0_{\tau_{t-}})^2) + \tilde \E^0((\eta^0_{\tau_{t+}} - \eta^0_{\tau_{t-}})^2) \leq 2 \;\tilde \E^0((2 \psi_t)^2) \label{firstterm}
\end{eqnarray}
and similarly,
\be{secondterm} \tilde \E^0((\eta^0_{\tau_{s+}})^2 - (\eta^0_{s})^2) \leq 2\;\tilde\E^0((2 \psi_s)^2). \ee
Finally, using (\ref{boundtaut}) and the convention $\tau_{r+} = 0$ when $r < 0$,
\begin{eqnarray}
\tilde \E^0((\eta^0_{\tau_{t+}}- \eta^0_{\tau_{s+}})^2) &=& \int_s^\infty \tilde \E^0((\eta^0_{\tau_{(t-r)+}})^2)\;\tilde \P^0(\tau_{s+} \in dr) \nonumber \\
&\leq& C\int_s^\infty ((t-r)\vee 0) \;\tilde \P^0(\tau_{s+} \in dr) \leq C(t-s). \label{thirdterm}
\end{eqnarray}
Using (\ref{firstterm}), (\ref{secondterm}) and (\ref{thirdterm}) back in (\ref{threeterms}), we are done.

$(iii.)$ For $l \geq 0$,
$$\P(|\eta^0_t| > l) = \P(t < \hat T^0 < \infty, |\eta^0_t| > l) + \P(\hat T^0 = \infty) \; \tilde \P^0(|\eta^0_t| > l).$$
The first term is less than
$$\P(\hat T^0 < \infty, M^0_{\hat T^0} > l) \leq \P(l/\kappa < \hat T^0 < \infty) + \P(M^0_{l/\kappa} > l),$$
where $\kappa$ is as in (\ref{comparePoi}). Now use (\ref{comparePoi}) and (\ref{fastDeath}) to get that this last sum is less than $Ce^{-cl}$. Next, we have
$$\tilde \P^0(|\eta^0_t| > l) \leq \tilde \P^0\left(\max_{1 \leq i \leq [t]} |\eta^0_{\tau^0_i}| > l/2\right) + \tilde \P^0(\psi_t > l/2),$$
again because there are at most $[t]$ renewals until time $t$. By (\ref{boundpsitail}), $\tilde \P^0(\psi_t > l/2) \leq Ce^{-cl}$. By Proposition 2.1.2 in \cite{lawler}, $\displaystyle{ \tilde \P^0\left(\max_{1 \leq i \leq [t]} |\eta^0_{\tau^0_i}| > l/2\right) \leq Ce^{-cl^2/[t]} }$. This completes the proof.
\epr

%%%%%%%%%%%%%% section 2 %%%%%%%%%%%%%%%%%%%%%%%%%%%%%%%%%%%%%%%%%%%%%%%%%%%%%%%%%%%%%%%%%%%%%%%%%%%%%%%%%%%%%%%%%%%%%%
\section{Pairs and sets of ancestries}

In this section, we study the joint behavior of ancestral paths. For pairs of ancestries, we define joint renewal points that have properties similar to the ones just discussed for single renewals, and then use these properties to study the speed of coalescence of first ancestrals. For sets of ancestries, we show that, given $N > 0$, the overall density of sites of $\Z$ occupied by ancestrals of rank smaller than or equal to $N$ at time $t$ tends to $0$ as $t \to \infty$.

Given $x, y \in \Z$, define $\tilde \P^{x, y}(\cdot) = \P(\cdot | \hat T^x = \hat T^y = \infty)$ and $\tilde \E^{x, y}(\cdot) = \E(\cdot |\hat T^x = \hat T^y = \infty)$. On $\{\hat T^x = \hat T^y = \infty\}$, let us define our sequence of \textit{joint renewal times}; start with $\tau^{x, y}_0 \equiv 0$,
$$\tau^{x,y}_1 = \inf\{t \geq 1: \hat T^{(\eta^x_t, t)} = \hat T^{(\eta^y_t, t)} = \infty\},$$
and, for $n \geq 1$, on $\{\hat T^x = \hat T^y = \infty, \eta^x_{\tau^{x, y}_n} = z, \eta^y_{\tau_n^{x,y}} = w\}$, define
$$\tau^{x,y}_{n+1} = \tau^{x, y}_n + \tau^{z, w}_1 \circ \theta(0, \tau^{x, y}_n).$$
Note that $\tilde \P^{x, x} = \tilde \P^x, \tilde \E^{x,x} = \tilde \E^x$ and $\tau^{x, x}_n = \tau^x_n$ for any $x$ and $n$. We have the following analog of Lemma \ref{renewal}:
\bp{jrenewal} $(i.)\; \tilde \P^{x, y}(\tau^{x, y}_n < \infty) = 1\; \forall n, x, y.$ \medskip\\
$(ii.)\;$ For $n \geq 0$, let
$$H_n = H_{[0, \tau^{x, y}_n(H)]}, \qquad H_{n+} = \theta(0, \tau^{x,y}_n)(H).$$
Given an event $A$ on finite-time Harris constructions, an event $B$ on Harris constructions and $z, w \in \Z$, we have
$$\tilde \P^{x,y}(H_n \in A, \eta^x_{\tau^{x, y}_n} = z,  \eta^y_{\tau^{x, y}_n} = w, H \in B) = \tilde \P^{x, y}(H_n \in A, \eta^x_{\tau^{x, y}_n} = z,  \eta^y_{\tau^{x, y}_n} = w) \cdot \tilde \P^{z, w}(H \in B).$$
$(iii.)\;$ Under $\tilde \P^{x, y}$, the $\Z^2$-valued process $\big(\eta^x_{\tau^{x,y}_n}, \eta^y_{\tau^{x,y}_n}\big)_{n \geq 0}$ is a Markov chain starting at $(x, y)$ and with transitions
$$P((a, b), (c, d)) = \tilde \P^{a, b}\big(\;\eta^a_{\tau^{a,b}_1} = c,\; \eta^b_{\tau^{a,b}_1 } = d \;\big).$$
In particular, if $\{\hat T^x = \hat T^y = \infty\}$ and $\eta^x_{\tau^{x,y}_m} = \eta^y_{\tau^{x,y}_m}$, then $\eta^x_{\tau^{x,y}_n} = \eta^y_{\tau^{x,y}_n}$ for all $n \geq m$.\\
$(iv.)\;$ There exist $c, C > 0$ such that, for any $x, y$,
$$\tilde \P^{x, y} \big(\max\big(\tau^{x, y}_1, M^x_{\tau^{x, y}_1}, M^y_{\tau^{x, y}_1}\big) > r\big) \leq Ce^{-cr}.$$
\ep
We omit the proof since it is an almost exact repetition of the one of Lemma \ref{renewal}; the only difference is that, when looking for renewals, we must inspect two points instead of one.

We now study the behavior of the discrete time Markov chain mentioned in part $(iii.)$ of the above proposition. Our first objective is to show that the time it takes for two ancestries to coalesce has a tail that is similar to that of the time it takes for two independent simple random walks on $\Z$ to meet. This fact will be extended to continuous time in Lemma \ref{aboutJ}; in Section 5, we will establish other similarities between pairs of ancestries and pairs of coalescing random walks.

\bl{discretecoal}
$(i.)$ For $z \in \Z$, let $\pi_z$ denote the probability on $\Z$ given by $$\pi_z(w) = \tilde \P^{0, z}\left(\eta^z_{\tau^{0, z}_1} - \eta^0_{\tau^{0, z}_1} = z + w\right), \quad w \in \Z.$$ There exist a symmetric probability $\pi$ on $\Z$ and $c, C > 0$ such that
$$||\pi_z - \pi||_{TV} \leq Ce^{-c|z|} \quad \forall z \in \Z,$$
where $||\cdot||_{TV}$ denotes total variation distance.\\
$(ii.)$ There exists $C > 0$ such that, for all $x, y \in \Z$ and $n \in \N$,
$$\tilde \P^{x,y}\big(\;\eta^x_{\tau^{x, y}_n} \neq \eta^y_{\tau^{x, y}_n} \;\big) \leq \frac{C|x-y|}{\sqrt{n}}.$$
\el
\bpr
For $n \geq 0, x, y \in \Z$, define
$$X^{x,y}_n = \left\{\begin{array}{ll} \eta^{y}_{\tau^{x,y}_n} - \eta^{x}_{\tau^{x,y}_n}, & \text{if } \hat T^x = \hat T^y = \infty; \\ \triangle, & \text{ otherwise}.\end{array} \right.$$
Using Proposition \ref{jrenewal}$(iii.)$ and translation invariance, we see that under $\tilde \P^{x, y}, X^{x,y}_n$ is a Markov chain that starts at $y-x$ and has transitions
$$\tilde \P^{x,y}(X^{x,y}_{n+1} = z + w \;|\; X^{x, y}_n = z) = \tilde \P^{0, z}(X^{0, z}_1 = z+ w) = \pi_z(w).$$
In particular, $0$ is an absorbing state. 

Fix $z \in \Z$ and $\kappa$ as in (\ref{comparePoi}). Let $t^* = z/3\kappa$. Let us take two random Harris constructions $H^1$ and $H^2$ defined on a common space with probability measure $\P$, under which $H^1$ and $H^2$ are independent and both have the original, unconditioned distribution obtained from the construction with Poisson processes. Define $H^3$ as a superposition of $H^1$ and $H^2$, as follows: we include in $H^3$:\\
$\bullet$ from $H^1$, all death marks in sites that belong to $(-\infty, [z/2]]$ and all arrows whose starting points belong to $(-\infty, [z/2]]$;\\
$\bullet$ from $H^2$, all death marks in sites that belong to $([z/2], \infty)$ and all arrows whose starting points belong to $([z/2], \infty)$.\\
Then, $H^3$ has same law as $H^1$ and $H^2$. We will write all processes and times defined so far as functions of these Harris constructions; for example, we will write $\eta^0_t(H^1)$ and consider $\tau^{0, z}_n(H^3)$ on the event $\{\hat T^0(H^3) = \hat T^z(H^3) = \infty\}$. Additionally, on the event $\{\hat T^0(H_1) = \hat T^z(H^2) = \infty\}$, define
$$\sigma^{0,z} = \inf\{t \geq 1: \hat T^{(\eta^0_t(H^1), t)}(H_1) = \hat T^{(\eta^z_t(H^2), t)}(H_2) = \infty\}.$$
As in Proposition \ref{jrenewal}$(iv.)$, there exist $c, C > 0$ such that
\be{speedpi}\P\big(\;\hat T^0(H^1) = \hat T^z(H^2) = \infty, \; \sigma^{0, z} \vee M^0_{\sigma^{0,z}}(H^1) \vee M^z_{\sigma^{0,z}}(H^2) > r \;\big) \leq Ce^{-cr}.\ee
Also define
$$Y^{0, z} =\left\{ \begin{array}{ll} \eta^z_{\sigma^{0, z}}(H^2) - \eta^0_{\sigma^{0, z}}(H^1), &\text{if } \hat T^0(H^1) = \hat T^z(H^2) = \infty, \\ \triangle, & \text{otherwise}.\end{array} \right.$$

Consider the events
$$\mathcal{L}_1 = \{M^0_{t^*}(H^1) \vee M^z_{t^*}(H^2) < z/2\},$$
$$\mathcal{L}_2 = \{\hat T^0(H^1) \wedge \hat T^z(H^2) < t^*\},$$
$$\mathcal{L}_3 = \left\{\begin{array}{c}\hat T^0(H^1) = \hat T^z(H^2) = \hat T^0(H^3) = \hat T^z(H^3) = \infty, \\ \tau^{0, z}_1(H^3) < t^*, \sigma^{0, z} < t^* \end{array}\right\}.$$
We claim that, if the event $\mathcal{L} := \mathcal{L}_1 \cap (\mathcal{L}_2 \cup \mathcal{L}_3)$ occurs, then $X^{0, z}_1(H^3) = Y^{0, z}$. To see this, assume first that $\mathcal{L}_1 \cap \mathcal{L}_2$ occurs. Then, we either have $\hat T^0(H_1) = \hat T^0(H_3) < t^* < \infty$ or $\hat T^z(H_2) = \hat T^z(H_3) < t^* < \infty$, and in either case $X^{0, z}_1(H^3) = Y^{0, z} = \triangle$. Now assume $\mathcal{L}_1 \cap \mathcal{L}_3$ occurs. Define
$$\begin{array}{ccc}t_1 = \tau_1^{0, z}(H^3), & a_1 = \eta^0_{t_1}(H^3), & b_1 = \eta^z_{t_1}(H^3), \\t_2 = \sigma^{0, z}, & a_2 = \eta^0_{t_2}(H^1), & b_2 = \eta^z_{t_2}(H^2). \end{array}$$
In $\mathcal{L}_1$, the ancestries of $0$ according to $H^1$ and $H^3$ coincide up to time $t^*$, and similarly for the ancestries of $z$ according to $H^2$ and $H^3$. Then, if we show that $t_1 = t_2$, we get $a_1 = a_2$ and $b_1 = b_2$, hence $X^{0, z}_1(H^3) = b_1 - a_1 = b_2 - a_2 = Y^{0, z}$. Assume $t_1 \leq t_2$. Since $\hat T^{(a_1, t_1)}(H^3) = \hat T^{(b_1, t_1)}(H^3) = \infty$, we have $\hat T^{(a_1, t_1)}(H^1), \hat T^{(b_1, t_1)}(H^2) > t^*$, so $a_2 = \eta^{(a_1, t_1)}_{1, t_2}(H^1), b_2 = \eta^{(b_1, t_1)}_{1, t_2}(H^2)$. But we also have $\hat T^{(a_2, t_2)}(H^1) = \hat T^{(b_2, t_2)}(H^2) = \infty$, so we get $\hat T^{(a_1, t_1)}(H^1) = \hat T^{(b_1, t_1)}(H^2) = \infty$, hence $t_1 = t_2$. By the same argument, $t_2 \leq t_1$ implies $t_2 = t_1$. This completes the proof of the claim.

Now note that the event $\mathcal{L}^c$ is contained in the union of:
$$\{M^0_{t^*}(H^1) > z/2\}, \{M^z_{t^*}(H^2) > z/2\},$$
$$\{t^* <\hat T^0(H^1) < \infty\}, \{t^* <\hat T^0(H^3) < \infty\}, \{t^* <\hat T^z(H^2) < \infty\}, \{t^* <\hat T^z(H^3) < \infty\},$$
$$\{\hat T^0(H^1) = \hat T^z(H^2) = \infty, \sigma^{0, z} > t^*\}, \{\hat T^0(H^3) = \hat T^z(H^3) = \infty, \tau^{0, z}_1(H_3) > t^*\}.$$
Using our choice of $t^*$, (\ref{comparePoi}), (\ref{fastDeath}), Proposition \ref{jrenewal}$(iii.)$ and (\ref{speedpi}), the probability of any of these events decreases exponentially with $z$.

We thus have
$\P(X^{0, z}_1(H^3) \neq Y^{0, z}) \leq Ce^{-cz}$, so
\be{coupling}\sum_{w \in \Z \cup \{\triangle\}} |\P(X^{0, z}_1(H^3) = w) - \P(Y^{0, z} = w)| \leq Ce^{-cz}.\ee
Now note that $\pi_z(\cdot) = \P(X^{0, z}_1(H^3) = z + \cdot \;|\; X^{0, z}_1(H^3) \neq \triangle)$ and define $\pi(\cdot) = \P(Y^{0, z} = z + \cdot \;|\; Y^{0, z} \neq \triangle)$. By the definition of $Y^{0, z}$ from independent Harris constructions, $\pi$ is symmetric and does not depend on $z$. Now,

$$\begin{aligned}&\sum_{w \in \Z} |\pi_z(w) - \pi(w)| = \sum_{w \in \Z} \left|\frac{\P(X^{0, z}_1(H^3) = w)}{\P(X^{0, z}_1(H^3) \neq \triangle)} - \frac{\P(Y^{0, z} = w)}{\P(Y^{0, z} \neq \triangle)}\right|\\
&\leq \frac{1}{\P(X^{0, z}_1(H^3) \neq \triangle)} \sum_{w \in \Z} |\P(X^{0, z}_1(H^3) = w) - \P(Y^{0, z} = w)| + \\&\qquad\qquad\qquad\left|\frac{1}{\P(Y^{0, z} \neq \triangle)} -  \frac{1}{\P(X^{0, z}_1(H^3) \neq \triangle)} \right| \sum_{w \in \Z}\P(Y^{0, z} = w).\\
&\leq \frac{1}{\P(X^{0, z}_1(H^3) \neq \triangle)} \; Ce^{-cz} + \frac{1}{\P(X^{0, z}_1(H^3) \neq \triangle)} \; \frac{1}{\P(Y^{0, z} \neq \triangle)} \; Ce^{-cz}.
\end{aligned}$$
We have $\P(X^{0, z}_1(H^3) \neq \triangle) = \P(\hat T^0(H^3) = \hat T^z(H^3) = \infty)$ and $\P(Y^{0, z} \neq \triangle) = \P(\hat T^0(H^1) = \hat T^z(H^2) = \infty)$, and these probabilities are bounded away from zero uniformly in $z$ by (\ref{poscorr}). We thus get $||\pi_z - \pi||_{TV} \leq Ce^{-cz}$ and part $(i)$ is proved.

We postpone the proof of part $(ii)$ -- the fact that for the perturbed random walk with increments $\pi_z$, the law of the hitting time of zero has same same tail as the one corresponding to a ``real'' random walk with increments $\pi$ -- to Section 6. Here, let us ensure that the four conditions in the beginning of that section are satisfied by $\pi_z$ and $\pi$. Conditions (\ref{symm}) and (\ref{totvar}) are already established. Condition (\ref{support}) is straightforward to check and (\ref{pitail}) follows from (\ref{speedpi}) and Proposition \ref{jrenewal}$(iv.)$.
\epr

We now want to define a random time $J^{x, y}$ that will work as a ``first renewal after coalescence'' for the first ancestrals of $x$ and $y$, a time after which the two processes evolve together with the law of a single first ancestor process. Some care should be taken, however, to treat the cases in which the ancestries of $x$ or of $y$ die out. With this in mind, we put
$$J^{x, y}= \left\{\begin{array}{ll}
\inf\{\tau^{x, y}_n: \eta^x_{\tau^{x,y}_n} = \eta^y_{\tau^{x,y}_n}\} &\text{on } \{\hat T^x = \hat T^y = \infty\};\medskip
\\ \inf\{\tau^x_n: \tau^x_n > \hat T^y\} &\text{on } \{\hat T^x = \infty, \hat T^y < \infty\}; \medskip
\\ \inf\{\tau^y_n: \tau^y_n > \hat T^x\} &\text{on } \{\hat T^y = \infty, \hat T^x < \infty\};\medskip
\\ 0 &\text{on } \{\hat T^x < \infty \text{ and } \hat T^y < \infty\}.\end{array} \right.$$
This definition is symmetric: $J^{x,y} = J^{y, x}$.
\bl{aboutJ}
$(i.)\;$ There exists $C > 0$ such that, for any $x, y \in \Z$ and $t \geq 0$,
$$\;\P(\;J^{x,y} > t \;) \leq \frac{C|x-y|}{\sqrt{t}}.$$
$(ii.)$ Conditioned to $\{\hat T^x=\infty\}$, the process $t \mapsto (\eta^0_{1, t}, \eta^0_{2,t}, \ldots) \circ \theta(\eta^x_{J^{x, y}}, J^{x,y})$ is independent of $J^{x,y}$. Additionally, the law of $t \mapsto (\eta^0_{1, t}, \eta^0_{2,t}, \ldots) \circ \theta(\eta^x_{J^{x, y}}, J^{x,y})$ conditioned to $\{\hat T^x = \infty\}$ is equal to the law of $t \mapsto (\eta^0_{1, t}, \eta^0_{2,t}, \ldots)$ conditioned to $\{\hat T^0 = \infty\}$.
\el
\bpr
$(i.)$ By Proposition \ref{jrenewal}$(iv.), \beta := 2 \sup_{z, w}\tilde \E^{z, w}(\tau^{z,w}_1)$ is finite. Noting that 
$$\tilde \P^{x,y}\left(\tau^{x,y}_{n+1}- \tau^{x,y}_n \in \cdot \;|\; \eta^x_{\tau^{x,y}_n} = z, \eta^y_{\tau^{x,y}_n} = w\right) = \tilde \P^{z,w}\left(\tau^{z, w}_1 \in \cdot\right),$$
Chebyshev's inequality yields $\tilde \P^{x,y}\left(\tau^{x,y}_{[t/\beta]} > t \right) \leq Ce^{-ct}$ for some $c, C > 0$, uniformly in $x, y$.
This together with Lemma \ref{discretecoal} gives
\be{fasttau*}\tilde \P^{x,y}(J^{x,y} > t) \leq \tilde \P^{x,y}\left(\eta^x_{\tau^{x, y}_{[t/\beta]}} \neq \eta^y_{\tau^{x, y}_{[t/\beta]}}\right) + \tilde \P^{x,y}\left(\tau^{x,y}_{[t/\beta]} > t\right) \leq \frac{C|x-y|}{\sqrt{t}}. \ee

Note that if $\hat T^x = \infty, \hat T^y < t/2$ and there exists some $n$ such that $\tau^x_n \in [t/2, t]$, then $J^{x, y} \leq t$, and similarly exchanging the roles of $x$ and $y$. Using (\ref{fastDeath}), Lemma \ref{exp0}$(i.)$ and (\ref{fasttau*}), we thus have
$$\begin{aligned}
&\P(J^{x,y} > t) \leq \\
&\P(\hat T^x = \infty, \hat T^y < \infty, J^{x, y} > t) + \P(\hat T^x < \infty, \hat T^y = \infty, J^{x, y}> t) + \P(\hat T^x = \hat T^y = \infty, J^{x, y} > t) \leq\\
&\P(t/2 < \hat T^y < \infty) + \P(\hat T^x = \infty, \nexists n: \tau^x_n \in [t/2, t]) +\\&\P(t/2 < \hat T^x < \infty) + \P(\hat T^y = \infty, \nexists n: \tau^y_n \in [t/2, t]) +\\&\qquad \P(\hat T^x = \hat T^y = \infty) \cdot \tilde \P^{x, y} (J^{x, y} > t) \leq \frac{C|x-y|}{\sqrt{t}}.
\end{aligned}$$

$(ii.)$ Let $A$ be a borelian of $[0, \infty)$ and $B$ be an event on Harris constructions. Using Proposition \ref{renewal} $(ii.)$,
$$\begin{aligned}
&\tilde \P^x(J^{x, y} \in A,\; \hat T^y < \infty,\; \theta(\eta^x_{J^{x, y}}, J^{x, y})(H) \in B) =\\
&\quad \sum_{n = 1}^\infty \tilde \P^x(\tau^x_{n-1} < \hat T^y \leq \tau^x_n,\; \tau^x_n \in A,\; \theta(\eta^x_{\tau^x_n}, \tau^x_n)(H) \in B) =\\
&\quad \tilde \P^0(H \in B) \cdot \sum_{n=1}^\infty \tilde \P^x(\tau^x_{n-1} < \hat T^y \leq \tau^x_n, \tau^x_n \in A) = \tilde \P^0(H \in B) \cdot \tilde \P^x(J^{x,y} \in A, \; \hat T^y < \infty)
\end{aligned}$$
Using Proposition \ref{jrenewal} $(ii)$ and the fact that $\tilde \P^{z, z} = \tilde \P^z$ for any $z$,
$$\begin{aligned}
&\tilde \P^x(J^{x, y} \in A, \; \hat T^y = \infty, \; \theta(\eta^x_{J^{x, y}}, J^{x, y})(H) \in B) =\\
&\frac{\P(\hat T^x = \hat T^y = \infty)}{\P(\hat T^x = \infty)} \sum_{n=1}^\infty \sum_{z \in \Z} \tilde \P^{x,y}\left(
\begin{array}{l}
\eta^x_{\tau^{x,y}_{n-1}} \neq \eta^y_{\tau^{x,y}_{n-1}}, \; \eta^x_{\tau^{x,y}_n} = \eta^y_{\tau^{x,y}_n} = z,\medskip\\ \; \tau^{x,y}_n \in A, \; \theta(z, \tau^{x,y}_n)(H) \in B 
\end{array}
\right) =\\
&\frac{\P(\hat T^x = \hat T^y = \infty)}{\P(\hat T^x = \infty)} \; \tilde \P^0(H \in B) \cdot \tilde \P^{x,y}(J^{x,y} \in A) = \tilde P^0(H \in B) \cdot \tilde \P^x(J^{x, y} \in A, \; \hat T^y = \infty).
\end{aligned}$$
Putting things together we get
$$\tilde \P^x(J^{x, y} \in A, \; \theta(\eta^x_{J^{x, y}}, J^{x, y})(H) \in B) = \tilde \P^0(H \in B) \cdot \tilde \P^x(J^{x, y} \in A).$$
The claim is a direct consequence of this equality.
\epr

\bl{rankCoal}There exist $c, C> 0$ such that, for any $x, y \in \Z, N \geq 1$ and $t \geq 0$,
$$\P(\hat T^x, \hat T^y > t, (\eta^x_{1, t}, \ldots, \eta^x_{N, t}) \neq (\eta^y_{1, t}, \ldots, \eta^y_{N, t})) \leq Ce^{CN - ct} + \frac{C|x-y|}{\sqrt{t}}.$$
\el
\bpr
There exists $\delta > 0$ such that, given a finite set $A \subset \Z$, we have $\P(\hat T^A < \infty) > \delta^{|A|}$. We can for instance take $\delta$ as the probability of a particle dying out before having any children, an observe that this occurs independently for different sites. This observation and the strong Markov property tell us that, defining
$\sigma_N = \sup\{s \geq 0: 0 < \#\eta^0_{*, s} < N\},$ we have $$\delta^N \P(\sigma_N > t) \leq \P(t < \hat T^0 < \infty).$$ Also using (\ref{fastDeath}), we obtain
\be{Nchildren}\P(\sigma_N > t) \leq C_1 e^{C_2 N - c_1 t}.\ee

Let $x, y \in \Z$; assume that $\hat T^x = \hat T^y = \infty$ and $J^{x,y} + \sigma_N \circ \theta(\eta^x_{J^{x,y}}, J^{x,y}) \leq t$. Noticing that $\eta^x_{J^{x, y}} = \eta^y_{J^{x,y}}$ and using the definition of $\sigma_N$, we must then have $z_1, \ldots, z_N$ such that $$\eta^0_{n, t - J^{x,y}} \circ \theta(\eta^x_{J^{x,y}}, J^{x,y}) = \eta^0_{n, t - J^{x,y}} \circ \theta(\eta^y_{J^{x,y}}, J^{x,y}) = z_n,\qquad 1 \leq n \leq N.$$ Lemma \ref{subdesc} then implies
$$\eta^x_{n, t} = \eta^y_{n, t} = z_n,\; 1 \leq n \leq N,$$
and thus
$$\begin{aligned} \P(\;\hat T^x &= \hat T^y = \infty,\;(\eta^x_{1, t}, \ldots, \eta^x_{N, t}) \neq (\eta^y_{1, t}, \ldots, \eta^y_{N, t})\;) \\
&\leq \P(\;\hat T^x = \hat T^y = \infty,\; J^{x,y} + \sigma_N \circ \theta(\eta^x_{1, J^{x,y}}, J^{x,y}) > t\;)\\
&\leq \P(\;J^{x,y} > t/2 \;) + \P(\;\sigma_N \circ \theta(\eta^x_{1, J^{x,y}}, J^{x,y}) > t/2\; \big|\; \hat T^x = \infty\;) \\
&\leq \frac{C|x-y|}{\sqrt{t}} + Ce^{CN - ct},\end{aligned}$$
where in the last inequality we used Lemma \ref{aboutJ}$(i.)$ in the first term and Lemma \ref{aboutJ}$(ii.)$ and (\ref{Nchildren}) in the second.

Finally, we have
$$\begin{aligned}
&\P(\hat T^x, \hat T^y > t, (\eta^x_{1, t}, \ldots, \eta^x_{N, t}) \neq (\eta^y_{1, t}, \ldots, \eta^y_{N, t}))\\
&\leq \P(t < \hat T^x < \infty) + \P(t < \hat T^y < \infty) + \P(\hat T^x = \hat T^y = \infty,(\eta^x_{1, t}, \ldots, \eta^x_{N, t}) \neq (\eta^y_{1, t}, \ldots, \eta^y_{N, t})) \\
&\leq 2Ce^{-ct} + Ce^{CN - ct} + \frac{C|x-y|}{\sqrt{t}} \leq Ce^{CN - ct} + \frac{C|x-y|}{\sqrt{t}}.\end{aligned}$$
\epr

\bp{densDecay}
There exist $C, \gamma > 0$ such that, for any $N \geq 1$ and $t \geq 0$,
$$\P(0 \in \{\eta^x_{n,t}: x \in \Z, 1 \leq n \leq N\}) \leq C \frac{N}{t^\gamma}.$$
\ep
\bpr
Fix a real $t \geq 0$ and a positive integer $l$ with $l > N$. Define $\Gamma = \{0, \ldots, l-1\}$ and $$\Lambda = \bigcap_{\{x,y\} \subset \Gamma} \left(\{\hat T^x \leq t\} \cup \{\hat T^y \leq t\} \cup \{(\eta^x_{1, t}, \ldots, \eta^x_{N,t}) = (\eta^y_{1, t}, \ldots, \eta^y_{N,t})\} \right).$$
We can use Lemma \ref{rankCoal} to bound the probability of $\Lambda^c$:
\be{LambdaComp}\P(\Lambda^c) \leq \sum_{\{x,y\} \subset \Gamma}\P(\hat T^x, \hat T^y > t, (\eta^x_{1, t}, \ldots, \eta^x_{N, t}) \neq (\eta^y_{1, t}, \ldots, \eta^y_{N, t}) ) \leq \frac{Cl^3}{\sqrt{t}} + Cl^2e^{CN - ct}\ee
since there are less than $l^2$ choices for $\{x, y\}$ and for any of them, $|x - y| \leq l$.

Let $\eta^\Gamma_{n,t} = \{\eta^x_{n,t} \in \Z: x \in \Gamma\}$ and $\eta^\Gamma_{N-, t} = \{\eta^x_{n, t} \in \Z: x \in \Gamma, 1 \leq n \leq N\}$. Since
$$\sum_{r=0}^{l-1}\sum_{n=1}^N \P(\Lambda, \{\eta^\Gamma_{n,t} \subset (r + l\Z)\}) \leq \sum_{n=1}^N \P(\Lambda) \leq N,$$
there exists $r^* \in \{0, \ldots, l-1\}$ such that 
\be{rstar}\sum_{n=1}^N \; \P\big(\;\Lambda, \; \{\eta^\Gamma_{n, t} \subset (r^* + l \Z)\}\;\big) \leq \frac{N}{l}. \ee

Finally, for $z \in \Z$ let $\Gamma_z = -r^* + lz + \Gamma$. The idea is that $0$ seen from $\Gamma_0$ is the same as $r^*$ seen from $\Gamma$. Let $\Lambda_z, \eta^{\Gamma_z}_{n, t}$ and $\eta^{\Gamma_z}_{N-, t}$ be defined from $\Gamma_z$ as $\Lambda, \eta^\Gamma_{n, t}$ and $\eta^\Gamma_{N-, t}$ are defined from $\Gamma$. Using (\ref{LambdaComp}) and (\ref{rstar}), we have

$$\begin{aligned}
&\P(0 \in \{\eta^x_{n, t}: x \in \Z, 1 \leq n \leq N\}) \leq \sum_{z \in \Z} \P(0 \in \eta^{\Gamma_z}_{N-, t}) \leq\\
&\sum_{z \in \Z} \P(0 \in \eta^{\Gamma_z}_{N-, t}, \Lambda_z) + \sum_{z \in \Z} \P(0 \in \eta^{\Gamma_z}_{N-, t}, \Lambda_z^c) \leq\\
&\sum_{z \in \Z} \sum_{n=1}^N \P(\Lambda_z, \; \{\eta^{\Gamma_z}_{n, t} = 0\} ) + \sum_{z \in \Z} \sum_{x \in \Gamma_z} \sum_{n=1}^N \P(\Lambda_z^c, \{\eta^x_{n, t} = 0\}) =\\
&\sum_{n=1}^N \sum_{z \in \Z} \P(\Lambda, \{\eta^\Gamma_{n, t} = r^* + lz\}) + \sum_{x \in \Gamma} \sum_{n=1}^N \sum_{z \in \Z} \P(\Lambda^c, \{\eta^x_{n, t} = r^* + lz\}) \leq \\
&\sum_{n=1}^N \P(\Lambda, \{\eta^\Gamma_{n, t} \subset (r^* + l \Z)\}) +\sum_{x \in \Gamma}^{l-1} \sum_{n=1}^N \P(\Lambda^c) \leq \frac{N}{l} + \frac{CNl^4}{\sqrt{t}} + CNl^3e^{CN - ct}.
\end{aligned}$$
We now put $l= t^{\frac{1}{9}}$ and observe that $(\frac{N}{t^{1/9}} + \frac{CNt^{4/9}}{t^{1/2}} + CNt^{1/3}e^{CN - ct}) \land 1 \leq \frac{CN}{t^\gamma}$ for some $C, \gamma > 0$.
\epr

%%%%%%%%%%%%%% section 4 %%%%%%%%%%%%%%%%%%%%%%%%%%%%%%%%%%%%%%%%%%%%%%%%%%%%%%%%%%%%%%%%%%%%%%%%%%%%%%%%%%%%%%%%%%%%%%

\section{Extinction and Survival}
In this section we prove Theorem \ref{extThm}. Our three ingredients will be a result about extinction under a stronger hypothesis (Lemma \ref{deathOfSingle}), an estimate for the edge speed of one of the types when obstructed by the other (Lemma \ref{speedObs}) and the formation of ``descendancy barriers'' for the contact process on $\Z$ (Lemma \ref{descBar}).

We recall our notation from the Introduction: the letters $\xi$ and $\eta$ will be used for the primal and dual process, respectively. Throughout this section, in contrast with the rest of the paper, Harris constructions and statements related to them, such as ``$(x, s) \leftrightarrow (y, t)$'', refer to the construction for the primal process rather than that of the dual.

\bl{deathOfSingle}
For the process $(\xi_t)$ with initial state $\xi_0$ such that $$\liminf_{x \rightarrow -\infty} \xi_0(x) = \liminf_{x \rightarrow \infty} \xi_0(x) = 2,$$ the $1$'s almost surely die out, \textit{i.e.} almost surely there exists $t$ such that $\xi_t(x) \neq 1 \; \forall x$.
\el
\bpr The hypothesis implies that there exists a finite $A \subset \Z$ such that $\xi_0(x) = 2 \; \forall x \in A^c$. Using Proposition \ref{densDecay}, we can choose $t$ such that $\mathcal{A} = \{A \cap \{\eta^x_{1, t}: x \in \Z\} \neq \emptyset\}$ has small probability. Place the primal time origin at dual time $t$; then, in $\mathcal{A}^c$, every site $x$ at primal time $t$ (\textit{i.e.} dual time $0$) either is in state $0$ or has its first ancestor $\eta^x_{1, t}$ in $A^c$, so $\xi_t(x) = 2$.
\epr

\bl{speedObs} Fix $\beta > 0$. For any $\epsilon > 0$, there exists $K > 0$ such that, if $\xi_0 = \xi^H = \mathds{1}_{(-\infty, 0]} + 2 \cdot \mathds{1}_{(0, \infty)}$, then
$$\P(\sup\{x: \xi_t(x) = 1\} \leq K + \beta t \; \forall t) > 1- \epsilon.$$
\el
\bpr
For $K > 0$, consider the events
\begin{eqnarray*}
A_n &=& \{\xi_n(x) = 1 \text{ for some } x \geq K/2 + \beta n/2\},\\
B_n &=& \{(x, n) \leftrightarrow (y, t) \text{ for some } x < K/2 + \beta n/2,\; y \geq K + \beta n,\; t \in [n, n+1]\},
\end{eqnarray*}
$n \in \{0, 1, 2, \ldots\}$. Now, using Lemma \ref{exp0} $(iii.)$,
$$\begin{aligned} 
\P(\cup_{n=0}^\infty A_n) &\leq \sum_{n=0}^\infty \; \sum_{x = K/2 + \beta n/2}^\infty \P(\eta^x_{1, n} \leq 0) \leq \sum_{n=0}^\infty \; \sum_{x = K/2 + \beta n/2}^\infty \P(|\eta^0_{1, n}| \geq x)\\
&\leq \sum_{n=0}^\infty \sum_{x = K/2 + \beta n/2}^\infty (Ce^{-cx^2/n} + Ce^{-cx}) \stackrel{K \to \infty}{\longrightarrow} 0.
\end{aligned}$$
Next, event $B_n$ requires the existence of a path that advances a distance of at least $K/2 + \beta n/2$ in a unit time interval; by a comparison with a sum of Poisson processes as in (\ref{comparePoi}), this occurs with probability smaller than $Ce^{-c(K/2 + \beta n/2)}$ for some $c, C > 0$, so $\P(\cup_n B_n) \leq \sum_n \P(B_n) \stackrel{K \to \infty}{\longrightarrow} 0$ as well. This gives $\P(\cap_n (A_n^c \cap B_n^c)) \to 1$ as $K \to \infty$, and to conclude the proof note that in $\cap_n (A_n^c \cap B_n^c)$, the set $\{(x, t): \xi_t(x) = 1\}$ is contained in $\{(x, t): x < K + \beta t\}$. 
\epr

For $\rho > 0$, define $V(\rho) = \{(x, t) \subset \Z \times [0, \infty): -\rho t \leq x \leq \rho t\}$. We say that site $0$ forms a $\rho$-descendancy barrier if\\
$(i.)$ for any $x, y \in \Z$ and $t \geq 0$ with $(x, 0) \leftrightarrow (y, t)$ and $(y, t) \in V(\rho)$, we have $(0, 0) \leftrightarrow (y, t)$;\\
$(ii.)$ for any $x, y \in \Z$ with opposite signs and $t \geq 0$ such that $(x, 0) \leftrightarrow (y, t)$, we have $(0, 0) \leftrightarrow (y, t)$.\\
Say that $x \in \Z$ forms a $\rho$-descendancy barrier if the origin forms a $\rho$-descendancy barrier according to $\theta(x, 0)(H)$.
\bl{descBar}
For any $\epsilon > 0$, there exists $\beta, K > 0$ such that
$$\P(\exists x \in [0, K]: x \text{ forms a } \beta \text{-descendancy barrier}) > 1 - \epsilon.$$
\el
The proof is in \cite{ampv}; see Proposition 2.7 and the definition of the event $\mathcal{H}_2$ in page 10 of that paper.

Finally, we state an obvious comparison result that can be verified by looking at the generator of the multitype contact process. As is usual, we abbreviate $\{x: \xi_t(x) = i\}$ as $\{\xi_t = i\}$.
\bl{monotone}
Let $(\xi '_t), (\xi ''_t)$ be two realizations of the multitype contact process built with the same Harris construction and such that $$\{\xi '_0 = 1\} \supset  \{\xi ''_0 = 1\}, \quad \{\xi '_0 = 2\} \subset \{\xi ''_0 = 2\}.$$ Then, $$\{\xi '_t = 1\} \supset  \{\xi ''_t = 1\}, \quad \{\xi '_t = 2\} \subset \{\xi ''_t = 2\} \quad  \forall t \geq 0.$$
\el
\bprthm{extThm}
We first prove that, if conditions $(i.)$ and $(ii.)$ in the statement of the theorem are satisfied, then the $1$'s become extinct. Fix $\epsilon > 0$. As in Lemma \ref{descBar}, choose $\beta, K_1$ corresponding to $\epsilon$, then as in Lemma \ref{speedObs}, choose $K_2$ corresponding to $\epsilon$ and $\beta$. Let $K = K_1 + K_2 + 2R$ (recall that $R$ is the range of the process). We may assume that there exist $a_1 < -L, a_2 > L$ (where $L$ is as in the statement of the theorem) such that $\xi_0(x) = 2 \; \forall x \in [a_1 - K, a_1] \cup [a_2, a_2 + K]$: after any positive time interval, there are infinitely many disjoint intervals of length $K$ that can be filled by a $2$ that is initially present. Let $(\xi^1_t), (\xi^2_t), (\xi^{12}_t)$ and $(\xi^{21}_t)$ be realizations of the multitype contact process all built using the same Harris construction as the original process $(\xi_t)$ and having initial configurations
$$\begin{aligned}
&\xi^1_0 = \mathds{1}_{(a_1, a_2)} + 2 \cdot \mathds{1}_{[a_1 - K, a_1]} + 2 \cdot \mathds{1}_{[a_2, a_2 + K]};\\
&\xi^2_0 = \mathds{1}_{(a_1, a_2)} + 2 \cdot \mathds{1}_{(a_1, a_2)^c};\\
&\xi^{12}_0 = \mathds{1}_{(-\infty, a_2)} + 2 \cdot \mathds{1}_{[a_2, \infty)};\\
&\xi^{21}_0 = 2 \cdot \mathds{1}_{(-\infty, a_1]} + \mathds{1}_{(a_1, \infty)}.
\end{aligned}$$
By a series of comparisons and uses of the previous lemmas, we will show that in $\xi^1$, the $1$'s become extinct with high probability. An application of Lemma \ref{monotone} to the pair $\xi^1, \xi$ then implies that in $\xi$, the $1$'s become extinct with high probability.

Define the events 
$$\mathcal{G}_1 = \{\forall t, \inf\{\xi^{21}_t= 1\} > a_1 - K_2 - \beta t\}, \quad
\mathcal{G}_2 = \{\forall t, \sup\{\xi^{12}_t = 1\} < a_2 + K_2 + \beta t\}\}.$$
By the choice of $K_2$, we have $\P(\mathcal{G}_1), \P(\mathcal{G}_2) > 1 - \epsilon$. Defining $W = \{(x, t): a_1 - K_2 - \beta t < x < a_2 + K_2 + \beta t\}$ and applying Lemma \ref{monotone} to the pairs $\xi^{12}, \xi^2$ and $\xi^{21}, \xi^2$, we get that
\be{stayInside}
\text{on } \mathcal{G}_1 \cap \mathcal{G}_2,\quad \{(x, t): \xi^2_t(x) = 1\} \subset W.
\ee
Also define
$$\mathcal{G}_3 = \{\exists b_1 \in [a_1 - K, a_1 - K + K_1]: b_1 \text{ forms a } \beta \text{-descendancy barrier}\};$$
$$\mathcal{G}_4 = \{\exists b_2 \in [a_2 + K - K_1, a_2 + K]: b_2 \text{ forms a } \beta \text{-descendancy barrier}\}.$$
The choice of $K_1$ and $\beta$ gives $\P(\mathcal{G}_3), \P(\mathcal{G}_4) > 1 - \epsilon$. Put $W_+ = \{(x, t): a_1 - K_2 - 2R - \beta t < x < a_2 + K_2 + 2R + \beta t\}$; a consequence of the definition of descendancy barriers is that
\be{barrier}
\text{on } \mathcal{G}_3 \cap \mathcal{G}_4,\quad \forall (x, t) \in W_+,\; \xi^1_t(x) = 0 \Leftrightarrow \xi^2_t(x) = 0.
\ee

We now claim that, in $\cap_{i=1}^4 \mathcal{G}_i, \{(x, t): \xi^1_t(x) = 1\} =\{(x, t): \xi^2_t(x) = 1\}$. This claim, together with Lemma \ref{deathOfSingle}, will imply that with probability larger than $1 - 4 \epsilon$, the $1$'s die out in $\xi^1$, and we will be done. To prove the claim, we start observing that $\{(x, t): \xi^1_t(x) = 1\} \supset \{(x, t): \xi^2_t(x) = 1\}$ always holds by Lemma \ref{monotone}. To establish the opposite inclusion in the occurrence of the good events, suppose to the contrary that for some $t, \{\xi^1_t = 1\} \neq \{\xi^2_t = 1\}$. But then we can find $(x^*, t^*)$ such that $\xi^1_{t^*}(x^*) = 1, \xi^2_{t^*}(x^*) \neq 1$ and $\{\xi^1_t = 1\} = \{\xi^2_t = 1\} \; \forall t \in [0, t^*)$. We must then have $\xi^1_{t^*-}(x^*) = 0$, since $\xi^1_{t^*-}(x^*) = 2$ would be incompatible with $\xi^1_{t^*}(x^*) = 1$ and $\xi^1_{t^*-}(x^*) = 1$ would imply, by the choice of $t^*, \xi^2_{t^*-}(x^*) = 1$ and then $\xi^2_{t^*}(x^*) = 1$, a contradiction. Now, since $\xi^1_{t^*-}(x^*) = 0$ and $\xi^1_{t^*}(x^*) = 1$ there must exist $y^*$ with $|y^* - x^*| \leq R$ such that $\xi^1_{t^*-}(y^*) = \xi^1_{t^*}(y^*) = 1$ and there exists an arrow from $(y^*, t^*)$ to $(x^*, t^*)$. But then, again by the choice of $t^*, \xi^1_{t^*-}(y^*) = 1$ implies $\xi^2_{t^*-}(y^*) = 1$, so $\xi^2_{t^*}(y^*) = 1$. Using (\ref{stayInside}), we can then conclude that $(y^*, t^*) \in W$, so $(x^*, t^*)$ is in the interior of $W_+$. This, (\ref{barrier}) and $\xi^1_{t^*-}(x^*) = 0$ imply that $\xi^2_{t^*-}(x^*) = 0$, so $\xi^2_{t^*}(x^*) = 1$, another contradiction. This completes the proof.

To prove the converse, we start noting that the case where there are infinitely many $1$'s in $\xi_0$ is trivial because then, at any $t \geq 0$ there almost surely exists some $x \in \Z$ such that $\xi_0(x) = 1$ and no death mark is present on $\{x\} \times [0, t]$, so the $1$'s are almost surely always present. We must thus show that, if condition $(i.)$ of the theorem is satisfied but condition $(ii.)$ is not, then the $1$'s have positive probability of surviving. By simple comparison arguments using Lemma $\ref{monotone}$, this reduces to proving that there exists $K > 0$ such that, if $\xi_0 = 2 \cdot \mathds{1}_{(-\infty, 0)} + \mathds{1}_{[0, K]}$, then $\P(\forall t, \{\xi_t = 1\} \neq \emptyset) > 0$. We will prove the stronger statement that this probability converges to $1$ as $K \to \infty$. Fix $\epsilon > 0$ and choose $\beta, K_1$ and $K_2$ as before. We will need another constant $K_3$ whose choice will depend on the following. Let $\alpha > 0$ be the edge speed for our contact process (i.e., the almost sure limit as $t \to \infty$ of $\frac{1}{t}\sup\{y: \exists x \in (-\infty, 0]:(x, 0) \leftrightarrow (y, t)\}$). Given $\alpha '  \in (0, \alpha)$, we have 
\be{alpha} \lim_{K' \to \infty} \P(\forall t, \;\exists x \in [0, K'], y > \alpha ' t: (x, 0) \leftrightarrow (y, t)) = 1.
\ee
This is a consequence of the definition of $\alpha$ and the fact that $\lim_{K' \to \infty} \P(\forall t, \; \exists x \in [0, K'], y \in \Z: (x, 0) \leftrightarrow (y, t)) = 1$; we omit the details. We may assume that the $\beta$ we have chosen is strictly smaller than $\alpha$, and we choose $K_3$ such that, putting $K' = K_3$ and $\alpha ' = \beta$, the probability in (\ref{alpha}) is larger than $1 - \epsilon$. Set $K = K_1 + K_2 + K_3 + 2R$. 

Recycling some of the notation from before, define $(\xi^{21}_t)$ with the same Harris construction as that of $(\xi_t)$, with
$$\xi^{21}_0 = 2 \cdot \mathds{1}_{(-\infty, 0)} + \mathds{1}_{[0, \infty)}$$
and the events
$$\begin{aligned}
&\mathcal{G}_1 = \{\forall t, \; \sup\{\xi^{21}_t = 2\} < K_2 + \beta t\};\\
&\mathcal{G}_2 = \{\exists x \in (K_2 + 2R, K_2 + 2R + K_1]: x \text{ forms a } \beta \text{-descendancy barrier}\};\\
&\mathcal{G}_3 = \{\forall t, \; \exists x \in (K_2 + 2R + K_1, K],\; y > K_2 + 2R + K_1 + \beta t: \; (x, 0) \leftrightarrow (y, t)\}.
\end{aligned}$$
We have $\P(\cap_{i=1}^3 \mathcal{G}_i) > 1 - 3\epsilon$. We can argue as before to the effect that, on $\mathcal{G}_1 \cap \mathcal{G}_2, \{\xi_t^{21} =2\} = \{\xi_t = 2\}$ holds for all $t$, so $\sup \{\xi_t = 2\} < K_2 + \beta t$ for all $t$. Additionally, on $\mathcal{G}_3$, for every $t$ there exists $y > K + \beta t$ such that $\xi_t(y) \neq 0$, so it must be the case that $\xi_t(y) = 1$. This shows that for all $t, \{\xi_t = 1\} \neq \emptyset$ and completes the proof.
\eprthm

%%%%%%%%%%%%%% section 5 %%%%%%%%%%%%%%%%%%%%%%%%%%%%%%%%%%%%%%%%%%%%%%%%%%%%%%%%%%%%%%%%%%%%%%%%%%%%%%%%%%%%%%%%%%%%%%
\section{Interface tightness}

We now carry out the proof outlined at the end of the Introduction. It is instructive to reestate Theorem \ref{interThm} in its dualized form:\\
\textbf{Theorem \ref{interThm}, dual version} For any $\epsilon > 0$, there exists $L > 0$ such that 
$$\P(|\sup\{x: \eta^x_t \leq 0\} - \inf\{x: \eta^x_t > 0\}| > L) < \epsilon \text{ for every } t \geq 0.$$

We start with two Lemmas concerning the expectation of the distance between two first ancestors. Lemma \ref{joinExp} shows a resemblance to the case of two random walks that evolve independently until they meet, at which time they coalesce. Lemma \ref{joinExpDeath} is a generalization that allows us to integrate over the event of death of a preassigned set of sites.

\bl{joinExp}There exists $C > 0$ such that, for all $x < y \in \Z$ and $t \geq 0$,\medskip\\
$(i.)\; \E\big(\; |\eta^y_t - \eta^x_t|\;) \leq C(y-x);$\medskip\\
$(ii.)\; \E\big(\;(\eta^y_t - \eta^x_t)^-) \leq C.$ 
\el
\bpr
By translation invariance, it suffices to treat $x = 0 < y$. It also suffices to prove $(i.)$ and $(ii.)$ for $t$ sufficiently large (not depending on $x, y$), because
$$\E\big(\;|\eta^y_t - \eta^0_t|\;\big) \leq y + \E\big(\;|\eta^0_t|\;\big) + \E\big(\;|\eta^y_t - y|\;\big) \leq y + \E(M^0_t) + \E(M^y_t) = y + 2 \E(M_t^0);$$
$$\E\big(\;(\eta^y_t - \eta^0_t)^-) \leq \E(\;(\eta^0_t)^+\;) + \E(\;(\eta^y_t - y)^-\;) \leq \E(M^0_t) + \E(M^y_t) = 2 \E(M_t^0),$$
and these expectations grow polynomially in $t$, by comparisons with sums of Poisson processes. Finally, 
\begin{eqnarray*}
\E \big(|\eta^0_t-\eta^y_t|\big) &=& \sum_{z, w} |z - w| \; \P\big(\eta^0_t = z, \; \eta^y_t = w\big)\\
&=& \sum_{z, w} |z - w| \; \P\big(\eta^0_t = z,\; \eta^y_t = w, \; \hat T^{(z, t)} = \hat T^{(w, t)} = \infty\big)\; \P\big(\hat T^{(z, t)} = \hat T^{(w, t)} = \infty\big)^{-1} \\
&\leq& C \; \E\big(|\eta^0_t - \eta^y_t|, \; \hat T^0 = \hat T^y = \infty\big)
\end{eqnarray*}
and similarly for $\E\big(\;(\eta^y_t - \eta^0_t)^-)$, so it suffices to prove $(i.)$ and $(ii.)$ on the event $\{\hat T^0 = \hat T^y = \infty\}$.\\
$(i.)$ We have
\begin{eqnarray}
&&\E\big(\;|\eta^{y}_t - \eta^{0}_t|; \; \hat T^0 = \hat T^y = \infty\;\big)   \nonumber \\&&\qquad\qquad\leq y + \E\big(\;|\eta^{0}_t|; \; \hat T^0 = \hat T^y = \infty,\; J^{0, y} > t\;\big) + \E\big(\;|\eta^{y}_t - y|; \; \hat T^0 = \hat T^y = \infty,\; J^{0, y} > t\;\big)\nonumber \\
&&\qquad\qquad= y + 2\E\big(\;|\eta^{0}_t|; \; \hat T^0 = \hat T^y = \infty,\; J^{0, y} > t\;\big) \label{joinExp2}
\end{eqnarray}
by symmetry. By Cauchy-Schwarz, this last expectation is less than
\be{joinExp3}
\left(\E\big((\eta^0_t)^2; \; \hat T^0 = \hat T^y = \infty,\; J^{0, y} > t \big) \cdot \P(\hat T^0 = \hat T^y = \infty,\; J^{0, y} > t)\right)^{\frac{1}{2}}.
\ee
Let us estimate the expectation.
$$\begin{aligned}&\E\big(\;(\eta^0_t)^2; \; \hat T^0 = \hat T^y = \infty,\; J^{0, y} > t \;\big) <\\
&\frac{1}{\P(\hat T^0 = \infty)} \cdot \E\big(\;(\eta^0_t)^2; \; \hat T^0 = \hat T^y = \infty, \; J^{0, y} > t\;\big) \leq \\
& \tilde \E^0\big(\;(\eta^0_t)^2; \; J^{0, y} > t\;\big)= \tilde \E^0((\eta^0_t)^2) - \tilde \E^0\big(\;(\eta^0_t)^2; \; J^{0, y} \leq t\;\big) =\end{aligned}$$
\be{joinExp11} \tilde \E^0\big((\eta^0_t)^2\big) - \tilde \E^0\big(\;(\eta^0_t - \eta^0_{J^{0,y}})^2 + (\eta^0_{J^{0,y}})^2 + 2 \eta^0_{J^{0, y}}(\eta^0_t - \eta^0_{J^{0,y}}); \; J^{0, y} \leq t \;\big).\ee
By Lemma \ref{aboutJ}$(ii.)$, we have
\begin{eqnarray}
&&\tilde \E^0\big(\;(\eta^0_t- \eta_{J^{0,y}})^2; \; J^{0, y} \leq t \;\big) = \int_0^t \tilde \E^0\big((\eta^0_{t-s})^2 \;\big) \cdot \tilde \P^0(J^{0, y} \in ds),\label{joinExp12} \\
&&\tilde \E^0\big(\;\eta^0_{J^{0,y}}(\eta^0_t- \eta^0_{J^{0,y}}); \; J^{0, y} \leq t \;\big) = 0.\label{joinExp13}
\end{eqnarray}
Using (\ref{joinExp12}) and (\ref{joinExp13}) and ignoring the term $(\eta^0_{J^{0, y}})^2$, the expression in (\ref{joinExp11}) is less than
\begin{eqnarray*}
&&\tilde \E^0\big((\eta^0_t)^2 \;\big) - \int_0^t \tilde \E^0\big((\eta^0_{t-s})^2 \;\big) \cdot \tilde \P^0(J^{0, y} \in ds) \\
&&\leq \tilde \E^0\big((\eta^0_t)^2\big) \cdot \tilde \P^0(J^{0, y} > t) + \int_0^t \tilde \E^0\big((\eta^0_t)^2 - (\eta^0_{t-s})^2\big) \cdot \tilde \P^0(J^{0, y} \in ds) \\
&&\leq (C_1 t + C_2) \frac{Cy}{\sqrt{t}} + \int_0^t (C_1 s + C_2) \; \tilde \P(J^{0, y} \in ds)
\end{eqnarray*}
by Lemma \ref{exp0}$(ii.)$ and Lemma \ref{aboutJ}$(i.)$. Now we can continue as in Lemma 1 in \cite{cd}: the above is less than
$$Cy \sqrt{t} + \frac{Cy}{\sqrt{t}} + C\int_0^t \tilde \P(J^{0, y} > u) \; du + C \leq Cy\sqrt{t} + C\int_0^t \frac{y}{\sqrt{u}}du \leq Cy\sqrt{t}$$
when $t \geq 1$. This and another application of Lemma \ref{aboutJ}$(i.)$ show that (\ref{joinExp3}) is less than $\sqrt{Cy\sqrt{t} \cdot \frac{Cy}{\sqrt{t}}} \leq Cy$; going back to (\ref{joinExp2}), we get
\be{joinExp4} \E\big(\;|\eta^{0}_t - \eta^{y}_t|; \; \hat T^0 = \hat T^y = \infty\big) \leq Cy. \nonumber \ee

$(ii.)$ To treat the expectation on the event $\{\hat T^0 = \hat T^y = \infty\}$, we will separately consider two cases, depending on whether or not the ancestor processes of $0$ and $y$ had a joint renewal in inverted order before time $t$. To this end, define
$$\tau^* = \inf\left\{\tau_n: \eta^{y}_{\tau_n^{0,y}} < \eta^{0}_{\tau_n^{0,y}}\right\}$$
(we set $\tau^* = \infty$ if the set is empty). Now,
\begin{eqnarray}
&&\E\big(\;(\eta^{y}_t - \eta^{0}_t)^-; \; \hat T^0 = \hat T^y = \infty, \; \tau^* \leq t \;\big) \nonumber \\
&&\leq \sum_{z < w} \int_0^t \tilde \E^{z, w}\big(\;|\eta^{w}_{t-s}-\eta^{z}_{t-s}| \;\big) \cdot \P\big(\;\hat T^0 = \hat T^y = \infty, \; \eta^{y}_{\tau^*} = z,\; \eta^{0}_{\tau^*} = w,\; \tau^* \in ds\;\big) \label{lastly} 
\end{eqnarray}
For each $z, w$, we have $\tilde \E^{z, w}(|\eta^w_{t-s} - \eta^z_{t-s}|) \leq \P(\hat T^z = \hat T^w = \infty)^{-1} \cdot C |w - z| \leq C |w - z|$ by part $(i.)$ and (\ref{poscorr}). Then, (\ref{lastly}) is less than
\begin{eqnarray}
&&C \sum_{z<w} \int_0^t(w-z)\; \P\big(\;\hat T^0 = \hat T^y = \infty,\; \eta^{y}_{\tau^*} = z,\; \eta^{0}_{\tau^*} = w\;; \tau^* \in ds\;\big) \nonumber \\
&&\leq C \sum_{z < w} (w-z) \; \P\big(\;\hat T^0 = \hat T^y = \infty, \; \eta^{y}_{\tau^*} = z,\; \eta^{0}_{\tau^*} = w,\; \tau^* < \infty\;\big)\nonumber \\ 
&&= C\; \E\big(\;(\eta^{y}_{\tau^*} - \eta^{0}_{\tau^*})^-;\;\hat T^0 = \hat T^y = \infty, \; \tau^* < \infty\;\big), \label{joinExp6}
\end{eqnarray}
which is bounded by Lemma \ref{discreteInv}.

Finally, as in Lemma \ref{exp0}, define on the event $\{\hat T^0 = \hat T^y = \infty\}$ the random variables $\tau^{0, y}_{t-}, \tau^{0, y}_{t+}$ and
$$\phi_t = M^{\left(\eta^{0}_{\tau^{0, y}_{t-}} ,\; \tau^{0, y}_{t-}\right)}_{\tau^{0, y}_{t+}} \vee M^{\left(\eta^{y}_{\tau^{0, y}_{t-}} ,\; \tau^{0, y}_{t-}\right)}_{\tau^{0, y}_{t+}}.$$
We then have $$\left|\eta^0_t - \eta^0_{\tau^{0, y}_{t-}}\right|, \left|\eta^y_t - \eta^y_{\tau^{0, y}_{t-}}\right| \leq \phi_t$$
on $\{\hat T^0 = \hat T^y = \infty\}$. Since on $\{\hat T^0 = \hat T^y = \infty, \tau^* > t\},\; \eta^{0}_{\tau^{0, y}_{t-}} \leq \eta^{y}_{\tau^{0, y}_{t-}}$ also holds, we have
\be{joinExp7} \E\big(\;(\eta^{y}_t - \eta^{0}_t)^-; \; \hat T^0 = \hat T^y = \infty, \; \tau^* > t\;\big) \leq  \E\big(\;2\phi_t; \; \hat T^0 = \hat T^y = \infty, \; \tau^* > t\;\big).\ee
As in the proof of Lemma \ref{exp0}, we can then show that $\E(\phi_t;\; \hat T^0 = \tilde T^y = \infty)$ is bounded uniformly in $y$ and $t$. Putting together (\ref{joinExp6}) and (\ref{joinExp7}), we get the result.
\epr

\bl{joinExpDeath}There exist $c, C > 0$ such that, for all $x < y \in \Z, t \geq 0$ and finite $A \subset \Z$,\medskip\\
$(i.)\; \E(\;|\eta^y_t - \eta^x_t|;\; \hat T^A < t\;) \leq C(y-x)e^{-c|A|};$\medskip\\
$(ii.)\; \E(\;(\eta^y_t - \eta^x_t)^-;\; \hat T^A < t\;) \leq Ce^{-c|A|}.$
\el
\bpr
Since both estimates are treated similarly, we will only show part $(ii.)$:
\begin{eqnarray}
&&\E(\;(\eta^y_t - \eta^x_t)^-;\; \hat T^A < t\;) \nonumber\\
&&\qquad = \sum_{k = 1}^\infty \E\left(\;(\eta^y_t - \eta^x_t)^-;\; \hat T^A < t, \; M^x_{\hat T^A} \vee M^y_{\hat T^A} = k\;\right) \nonumber \\
&& \qquad\leq \sum_{k=1}^\infty \sum_{i=-k}^k \sum_{j=-k}^k \E\left(\;\left(\eta^{(y+j, \hat T^A)}_{t - \hat T^A} - \eta^{(x+i, \hat T^A)}_{t - \hat T ^A}\right)^-;\; \hat T^A < t, \; M^x_{\hat T^A} \vee M^y_{\hat T^A} = k\;\right) \nonumber \\
&&\qquad \leq \sum_{k=1}^\infty \int_0^t \left(\sum_{i=-k}^k\sum_{j=-k}^k\E((\eta^{y+j}_{t-s} - \eta^{x+i}_{t-s})^-) \right)\P(\; \hat T^A \in ds, M^x_{\hat T^A} \vee M^y_{\hat T^A} = k\;) \nonumber
\end{eqnarray}
If $x+i < y +j$, then $\E\big(\;\big(\eta^{y+j}_{t-s} - \eta^{x+i}_{t-s}\big)^-\;\big) \leq C$ by Lemma \ref{joinExp}$(ii.)$. If $x+i > y+j$, then we also have $(x+i)-(y+j) < 2k$, so $\E\big(\;\left| \eta^{y+j}_{t-s} - \eta^{x+i}_{t-s} \right|\;\big) \leq 2Ck$ by Lemma \ref{joinExp}$(i.)$. Hence, in all cases the expectation is less than $Ck$, and the above sum is less than
\begin{eqnarray}
&& C \sum_{k=1}^\infty k^3 \; \P(\;\hat T^A < t,\; M^x_{\hat T^A} \vee M^y_{\hat T^A} = k\;) \leq C\; \E\left(\;(M^x_{\hat T^A} \vee M^y_{\hat T^A})^3; \; \hat T^A < \infty\right) \nonumber\\
&&\qquad \leq C\; \E\left((M^x_{\hat T^A})^3;\; \hat T^A < \infty\right) + C\; \E\left((M^y_{\hat T^A})^3;\; \hat T^A < \infty\right) \label{joinExpDeath1}.
\end{eqnarray}
Now, by Cauchy-Schwarz,
\be{joinExpDeath2}\E\left((M^x_{\hat T^A})^3;\; \hat T^A < \infty\right) \leq \left( \E\left((M^x_{\hat T^A})^6;\; \hat T^A < \infty\right)\cdot \P\big(\;\hat T^A < \infty\;\big) \right)^{1/2}.\ee
The probability in the right-hand side decreases exponentially with $|A|$ (see Section 11b in \cite{d1}). Doing
$$\P\big(\;M^x_{\hat T^A} > l, \; \hat T^A < \infty \;\big) \leq \P\left(\frac{l}{\sigma} < \hat T^A < \infty \right) + \P\left(\;M^0_{l/\sigma} > l\;\right)$$
with large $\sigma$ and using (\ref{fastDeath}) again, we see that the expectation on the right-hand side of (\ref{joinExpDeath2}) is uniformly bounded in $x$ and $A$. 
\epr

For $z > 0$, say that sites $x, x + z$ produce a $z$-inversion at time $t$ if $\eta^x_t > 0 \geq \eta^{x+z}_t$. The following lemma shows that the expected number of $z$-inversions at time $t$ is bounded uniformly in $z$ and $t$. It also illustrates the usefulness of Lemma \ref{joinExpDeath}.
\bl{addprob} There exist $c, C > 0$ such that, for any integer $z \geq 1$, real $t \geq 0$ and finite $A \subset \Z$,\medskip \\
$(i.)\; \displaystyle{\sum_{x \in \Z}} \; \P(\;\eta^x_{t} > 0 \geq \eta^{x+z}_{t},\; \hat T^{x+A} < t \;) \leq Ce^{-c|A|}$; \medskip\\
$(ii.)\; \displaystyle{\sum_{x \in \Z}} \;\P(\;\eta^x_{t} \leq 0 < \eta^{x+z}_{t},\; \hat T^{x+A} < t \;) \leq C|z|e^{-c|A|}.$
\el
\bpr
We start proceeding like in Lemma 4 in \cite{cd}, noticing that, by translation invariance,
$$\P(\; \eta^x_{t} > 0 \geq \eta^{x+z}_{t},\; \hat T^{x+A} < t \;) = \P(\; \eta^0_{t} > -x \geq \eta^{z}_{t},\; \hat T^{A} < t\;),$$
$$\P(\; \eta^x_{t} \leq 0 < \eta^{x+z}_{t},\; \hat T^{x+A} < t\;) = \P(\;\eta^0_{t} \leq -x < \eta^{z}_{t},\; \hat T^A < t \;)$$
and summing over $x$ to obtain
$$\sum_{x \in \Z}\P(\; \eta^x_{t} > 0 \geq \eta^{x+z}_{t},\; \hat T^{x+A} < t \;)  = \E\big(\;(\eta^z_{t} - \eta^0_{t})^-;\; \hat T^A < t \;),$$
$$\sum_{x \in \Z}\P(\; \eta^x_{t} \leq 0 < \eta^{x+z}_{t},\; \hat T^{x+A} < t \;) = \E\big(\;|\eta^0_{t} - \eta^z_{t}|;\; \hat T^A < t \big);$$
see Lemma 4 in \cite{cd} for more details. Also recall our conventions about the $\triangle$ state in Remark \ref{triangle}. Now, it suffices to apply Lemma \ref{joinExpDeath}.
\epr

Fix $0 < s < t$. For $x \in \Z$ such that $\eta^x_{*, t} \neq \emptyset$, let $n$ be the smallest integer such that $(\eta^x_{n, s}, s)$ survives up to time $t$ (as in the statement of Lemma \ref{subdesc}). Define $R^x(s,t) = \eta^x_{n, s}$.  Assume the primal time origin is at dual time $t$; since $\eta^x_{t} = \eta^{(R^x(s,t),s)}_{t}$, we have $\xi_t(x) = \xi_{t-s}(R^x(s, t))$.  

Also define $R(s,t) = \{R^x(s,t): x \in \Z, \eta^x_{*, t} \neq \emptyset\}$. This will be understood as a set of ``relevant'' sites. To get some insight into this, again assume that the primal time origin is placed at dual time $t$. Fix $y$ such that $\xi_{t-s}(y) \neq 0$ and change $\xi_{t-s}$ in the following way: switch the type of the individual at $y$ to the opposite one, and leave other sites untouched. Then let this new configuration evolve following the original primal Harris construction from primal time $t - s$ to $t$; denote by $\tilde \xi^y$ the final configuration obtained. Then, $R(s, t)$ is exactly the set of occupied sites $y$ in $\xi_{t-s}$ for which $\tilde \xi^y \neq \xi_t$. 

Our next task is to show that, if $s$ is large, then with high probability the restriction of $\xi_{t-s}$ to $R(s,t)$ has no interface. Formally, 

\bp{AKT} Let
\be{goodevent}\mathcal{G}(s, t) = \left\{\sup\{x \in R(s,t): \eta^{(x, s)}_t \leq 0\} < \inf\{x \in R(s, t): \eta^{(x, s)}_t > 0\}\right\}. \ee
Then, $\displaystyle \lim_{s \rightarrow \infty} \inf_{t \geq s} \P(\mathcal{G}(s,t)) = 1.$
\ep
\bpr
We fix $s < t$ and an integer $N$ to be chosen later. We will write $\mathcal{G}, R$ instead of $\mathcal{G}(s,t), R(s,t)$, and in general omit the dependence on $s, t, N$. 

Fix $d$ with  $1 > d \geq \P(0 \in \{\eta^x_{n,s}: x \in \Z, 1 \leq n \leq N\})$ and let $X$ be a random variable with uniform distribution on $\{0, \ldots, \lceil 1/d \rceil-1\}$ and independent of the Harris construction. Define
$$\hat R =  \{\eta^x_{n,s}: x \in \Z, 1 \leq n \leq N\} \cup (X + \lceil 1/d \rceil \Z).$$
$\hat R$ is a random subset of $\Z$; its law is invariant with respect to shifts in $\Z$ and $\P(0 \in \hat R) \leq 2d.$ Additionally, it only depends on the Harris construction on times in $[0, s]$, and of course on $X$. Put $\mathcal{S} = \{x \in \Z: \eta^{(x,s)}_{*, t} \neq \emptyset\}$. Note that by the definition of $R$, we have $R \subset \mathcal{S}$; also, by our conventions, when we say for example $\eta^{(x, s)}_{t} \geq 0$, we are implying that $x \in \mathcal{S}$.

We will also need the events
$$\mathcal{G}_1 = \{\;\nexists \; x, y \in \hat R: x < y,\; (x, y) \cap \hat R \cap \mathcal{S} = \emptyset,\; \eta^{(x,s)}_{t} > 0 \geq \eta^{(y, s)}_{t} \;\},$$
$$\mathcal{G}_2 = \{\;\nexists \; x \in R - \hat R,\; y \in \Z: x < y,\; (x, y) \cap \hat R\cap \mathcal{S} = \emptyset,\; \eta^{(x,s)}_{t} > 0 \geq \eta^{(y, s)}_{t}\;\},$$
$$\mathcal{G}_3 = \{\;\nexists \; y \in R - \hat R,\; x \in \Z: x < y,\; (x, y) \cap \hat R\cap \mathcal{S} = \emptyset,\; \eta^{(x,s)}_{t} > 0 \geq \eta^{(y, s)}_{t}\;\}.$$

We claim that $\mathcal{G}_1 \cap \mathcal{G}_2 \cap \mathcal{G}_3 \subset \mathcal{G}.$ Indeed, assume the three events occur and let us show that, given $a \in R$ such that $\eta^{(a, s)}_{t} > 0$, we have $\eta^{(b, s)}_{t} > 0$ for any $b > a, b \in R$. Let $\{z_1, z_2, \ldots\} = [a, \infty) \cap \hat R \cap \mathcal{S}$ with $z_i \leq z_{i+1} \; \forall i$. If $a < z_1$, then $a \in R - \hat R$, so $\eta^{(b, s)}_{t} > 0$ for any $b \in (a, z_1]$ by the definition of $\mathcal{G}_2$. If $a = z_1$, then we plainly have $\eta^{(z_1, s)}_{t} > 0$. So in any case we have $\eta^{(z_1, s)}_{t} > 0$, and from this we can use the definition of $\mathcal{G}_1$ to conclude that $\eta^{(z_i, s)}_{t} > 0 \; \forall i.$ Finally, if $b > z_1, b \in R$, then either $b = z_i$ for some $i$ or $b \in (z_i, z_{i+1})$ for some $i$. In the first case, we already have $\eta^{(b, s)}_{t} > 0$; in the second case, we have $b \in R - \hat R$, so we can apply the definition of $\mathcal{G}_3$ to $z_i$ and $b$ to conclude that $\eta^{(b, s)}_{t} > 0$. This concludes the proof of the claim.

Let us now estimate the probabilities of $\mathcal{G}_1^c, \mathcal{G}_2^c$ and $\mathcal{G}_3^c$.  
\begin{eqnarray*}
\P(\mathcal{G}_1^c) &\leq& \sum_{x < y} \P(\;[x, y] \cap \hat R \cap \mathcal{S} = \{x, y\},\; \eta^{(x,s)}_{t} > 0 \geq \eta^{(y, s)}_{t} \;)\\
&\leq& \sum_{x \in \Z, z \geq 1}\; \sum_{A \subset (0, z)}\P\big(\;[x, x+z] \cap \hat R = \{x, x+z\} \cup (x+A),\\
&&\qquad\qquad\qquad\qquad\qquad x + A \subset \mathcal{S}^c, \eta^{(x,s)}_{t} > 0 \geq \eta^{(x+z, s)}_{t}\;\big)\\
&=& \sum_{z, A} \P\big(\;[0, z] \cap \hat R = \{0, z\} \cup A \;\big) \cdot\\
&&\qquad\qquad\qquad\qquad\qquad \sum_{x\in \Z}\P\big(\;\hat T^{x+A} < t-s, \eta^{x}_{t-s} > 0 \geq \eta^{x+z}_{t-s}\;\big).
\end{eqnarray*}
Applying Lemma \ref{addprob} to the inner sum, we get that the above is less than
\begin{eqnarray*}
&C&\sum_{z, A} e^{-c(\#A)} \; \P\big(\;[0, z] \cap \hat R = \{0, z\} \cup A \;\big)\\
&\leq& C \sum_{k\geq 0} e^{-ck} \sum_{z \geq k+1}\; \sum_{A \subset (0, z): \#A=k}\P\big(\;[0, z] \cap \hat R = \{0, z\} \cup A \;\big)\\
&=& C \sum_{k \geq 0} e^{-ck} \;\P \big( \; 0 \in \hat R\; \big) \leq Cd.
\end{eqnarray*}

Similarly,
\begin{eqnarray}
\P(\mathcal{G}_2^c) &\leq& \sum_{x<y} \P\big(\;x \in R - \hat R,\; (x, y) \cap \hat R \cap \mathcal{S} = \emptyset,\; \eta^{(x,s)}_{t} > 0 \geq \eta^{(y, s)}_{ t} \;\big) \nonumber\\
&\leq& \sum_{x \in \Z, z \geq 1}\; \sum_{A \subset (0, z)}\; \sum_{a \in \Z , m > N} \;\sum_{(a_1, \ldots, a_{m-1}) \in \Z^{m-1}} \nonumber\\
&&\qquad\qquad\qquad\qquad\qquad \P\big(\;x = \eta^{x+a}_{m,s},\; x + a_i = \eta^{x+a}_{i,s} \; \forall i < m,\nonumber\\
&&\qquad\qquad\qquad\qquad\qquad (x, x+z) \cap \hat R = x + A,\nonumber\\
&&\qquad\qquad\qquad\qquad\qquad x + a_i \notin \mathcal{S} \; \forall i < m,\; x + A \subset \mathcal{S}^c,\; \eta^{(x,s)}_{t} > 0 \geq \eta^{(x+z, s)}_{t} \;\big)\nonumber\\
&\leq& \sum_{z, A, a, m, (a_i)} \P\big(\;0 = \eta^a_{m, s},\; a_i = \eta^a_{i, s} \; \forall i < m, \; (0, z) \cap \hat R = A\;\big) \cdot \nonumber\\
&&\qquad\qquad\quad \sum_{x \in \Z}\P\big(\;\hat T^{x+ a_i} < t-s \; \forall i < m, \; \hat T^{x+A} < t-s, \; \eta^{x}_{t-s} > 0 \geq \eta^{x+z}_{t-s}\;\big)\nonumber\\
&\leq& C\sum_{z, A, a, m, (a_i)}e^{-c ((\#A) \vee m)}\; \P\big(\;0 = \eta^a_{m, s},\; a_i = \eta^a_{i, s} \; \forall i < m, \; (0, z) \cap \hat R = A\;\big)\nonumber\\
&\leq& C\sum_{z, A, a, m}e^{-c ((\#A) \vee m)}\; \P\big(\;0 = \eta^a_{m, s},\; (0, z) \cap \hat R = A\;\big)\nonumber\\
&=& C \sum_{k \geq 0} \;\sum_{m > N}e^{-c (k \vee m)} \;\sum_{a \in \Z}\; \sum_{z \geq k+1} \; \sum_{A \subset (0, z): \#A = k} \P\big(\;0 = \eta^a_{m,s},\; (0, z) \cap \hat R = A\;\big)\nonumber\\
&=& C \sum_{k \geq 0} \;\sum_{m > N}e^{-c (k\vee m)} \;\sum_{a \in \Z}\; \sum_{z \geq k+1} \P\big(\;0 = \eta^a_{m,s},\; \#((0, z) \cap \hat R) = k\;\big). \label{huge}
\end{eqnarray}
Now note that, since $X + \lceil 1/d  \rceil \Z \subset \hat R$, there are no intervals of length larger than $\lceil 1/d \rceil$ that do not intersect $\hat R$. Hence, when $z > \frac{k+2}{d}$, we have $\#((0, z) \cap \hat R) > k$, hence $\P(\#((0, z) \cap \hat R) = k) = 0$. When $z \leq \frac{k+2}{d},$ we use the bound $\P(0 = \eta^a_{m,s}, \#((0, z) \cap \hat R) = k) \leq \P(0 = \eta^a_{m,s}).$ So the expression in (\ref{huge}) is less than
\be{lesshuge}C \sum_{k \geq 0} \;\sum_{m > N}\frac{k+2}{d}e^{-c (k \vee m)} \;\sum_{a \in \Z}\; \P\big(\;0 = \eta^a_{m,s}\;\big).\ee
The inner sum is less than
$$\sum_{a \in \Z} \P\big(\;0 \in \eta^a_{*,s} \;\big) = \E\;\#\{a \in \Z: 0 \in \eta^a_{*,s}\}.$$
By a routine comparison with Poisson process, the latter is less than $Cs$ for some $C > 0$. Hence the expression in (\ref{lesshuge}) is less than
$$\frac{Cs}{d} \sum_{k \geq 0} \sum_{m \geq N} (k+2) e^{-c(k \vee m)} \leq \frac{Cs}{d}\sum_{k \geq 0} (k+2)e^{-(c/2)k}\sum_{m \geq N}e^{-(c/2)m} \leq C\frac{s}{d}e^{-cN}$$
for some $c, C > 0$.

By symmetry, we have $\P\{\mathcal{G}_3\} =\P\{\mathcal{G}_2\}$. To summarize, we obtained:
\be{boundG1}\P(\mathcal{G}_1^c) \leq Cd;\ee
\be{boundG23}\P(\mathcal{G}_2^c), \P(\mathcal{G}_3^c) \leq C\frac{s}{d}e^{-cN}. \ee
Additionally, remember that we chose $d$ satisfying
\be{choiced}\P(0 \in \{\eta^{x}_{n, s}: x \in \Z, 1 \leq n \leq N\}) \leq d\ee
and Proposition \ref{densDecay} tells us that
\be{dDrecap} \P(0 \in \{\eta^{x}_{n, s}: x \in \Z, 1 \leq n \leq N\}) \leq C\frac{N}{s^\gamma}.\ee
So, putting $N = \lceil s^{\gamma/2} \rceil$ and $d = C\frac{N}{s^\gamma}$ (provided $s$ is large enough so that this is less than $1$), we conclude that $\P(\mathcal{G}^c) \leq \P(\mathcal{G}_1^c \cup \mathcal{G}_2^c \cup\mathcal{G}_3^c) \leq \P(\mathcal{G}_1^c)+ \P(\mathcal{G}_2^c)+ \P(\mathcal{G}_3^c) \to 0$ as $s \to \infty$.
\epr

Following the terminology in \cite{cd}, define $B_t = \#\{(x, y) : x < y, \eta^x_{t} > 0 \geq \eta^y_{t}\}$. Our next-to-last result before the proof of Theorem \ref{interThm} will be

\bp{expinv} The process $(B_t)_{t \geq 0}$ is tight.
\ep
\bpr
Let $\epsilon > 0$. By Proposition \ref{AKT}, there exists $s$ such that $\P(\mathcal{G}(s, t)^c) < \epsilon/2$ for any $t > s$. Fix $t > s$; we have
\begin{eqnarray}
\E(B_t; \mathcal{G}(s, t)) &=& \sum_{a < b} \P\big(\;\eta^{a}_{t} > 0 \geq \eta^{b}_{t},\; \mathcal{G}(s, t)\;\big) \nonumber\\
&\leq&\sum_{a < b}\;\sum_{x < y} \P\big(\;R^a_{s, t} = y, \; R^b_{s, t} = x,\; \eta^{(x, s)}_{t} \leq 0 < \eta^{(y, s)}_{t}\;\big)\nonumber\\
&\leq&\sum_{a < b}\;\sum_{x < y} \P\big(\;y \in \eta^{a}_{*, s}, \;  x \in \eta^{b}_{*,s},\; \eta^{(x, s)}_{t} \leq 0 < \eta^{(y, s)}_{t}\;\big)\nonumber\\
&=& \sum_{z \geq 1} \; \sum_{x \in \Z} \P\big(\;\eta^{(x, s)}_{t} \leq 0 < \eta^{(x+z, s)}_{t}\;\big)  \; \sum_{a < b} \P\big(\;z \in \eta^{a}_{*, s}, \;  0 \in \eta^{b}_{*,s}\;\big). \label{interBt}
\end{eqnarray}
By (\ref{comparePoi}), there exist $c$ (that depends on $s$) such that 
$$\P\big(\;z \in \eta^a_{*, s}\;\big) \land \P\big(\;0 \in \eta^b_{*, s}\;\big) \leq \P(M^0_s > |a-z|) \land \P(M^0_s > |b|) \leq e^{-c(|a-z| \vee |b|)},$$
then
$$\sum_{a < b} \P\big(\;z \in \eta^{a}_{*, s}, \;  0 \in \eta^{b}_{*,s}\;\big) \leq \sum_{a < b} (e^{-c(|a-z| \vee |b|)}) \leq Ce^{-cz}$$
as is easily seen. Using this and Lemma \ref{addprob}, we see that the expression in (\ref{interBt}) is less than
$$C\sum_{z \geq 1} e^{-cz}\; \sum_{x \in \Z} \P\big(\;\eta^{(x, s)}_{t} < 0 \leq \eta^{(x+z, s)}_{t}\;\big) \leq C \sum_{z \geq 1} ze^{-cz} < \infty.$$
So, if $L$ is large, we have $\frac{\E(B_t; \mathcal{G}(s, t))}{L} < \frac{\epsilon}{2}$ for all $t > s$, and thus
$$\P(B_t > L) \leq \P(\mathcal{G}(s, t)^c) + \P(B_t > L, \mathcal{G}(s, t)) \leq \frac{\epsilon}{2} + \frac{\E(B_t; \mathcal{G}(s, t))}{L} < \epsilon,$$
Noticing that the trajectories of $(B_t)$ are right continuous with left limits, we can increase $L$ if necessary so that this inequality also holds for $t \leq s$, completing the proof.
\epr

\bprthm{interThm}
We separately show that $(\rho_t \land 0)$ and $(\rho_t \vee 0)$ are tight. We start with the first. Given $L > 0$, for the event $\{\rho_t > L\}$ to occur, there necessarily exist two sites $x, y$ such that $y - x > L$ and $\eta^y_{t} \leq 0 < \eta^x_{t}$. If $N < L$ and $\{B_t < N\}$ also occurs, then we cannot have more than $N$ sites $z \in (x, y)$ such that $\eta^z_{*, t} \neq \emptyset$, because every such site produces a crossing either with $x$ or with $y$ and thus increases $B_t$ by one. So we have, for all $t \geq 0$,
$$\begin{aligned}
\P(B_t < N, \rho_t > L) &\leq \sum_{x < y, y - x > L}\P\big(\;\eta^x_{t} > 0 \geq \eta^y_{t},\; \hat T^{(x, y)\backslash A} < t \text{ for some } A \subset (x, y), \#A < N\;\big)\\
&\leq \sum_{z > L} \; \sum_{A \subset (0, z): \#A < N} \; \sum_{x \in \Z} \P\big(\;\eta^x_{t} > 0 \geq \eta^{x+z}_{t},\; \hat T^{(x, x + z) \backslash (x + A)} < \infty\big).
\end{aligned}$$
Using Lemma \ref{addprob} on the innermost sum and counting the possible choices of $A$, the above is less than
$$C\sum_{z > L} \left[{z \choose 0} + \ldots + {z \choose N-1} \right]e^{-c (z - N)},$$
which tends to $0$ as $L \to \infty$. So, given $\epsilon > 0$, choose $N > 0$ such that $\P(B_t \geq N) < \epsilon/2 \; \forall t$, then choose $L$ such that $\P(B_t < N, \rho_t > L) < \epsilon/2\; \forall t$, so that $\P(\rho_t > L) \leq \P(B_t \geq N) + \P(B_t < N, \rho_t > L) < \epsilon \; \forall t$, and we are done.

Now we treat $(\rho_t \vee 0).$ This is easier: given $L > 0$, for $\{\rho_t < -L\}$ to occur we must have $x < y$ such that $\eta^{x}_{t} < 0 \leq \eta^{y}_{t}$ and $\eta^{w}_{*, t} = \emptyset \; \forall w \in (x, y)$. Then, for any $t$,
$$\begin{aligned}
\P(\rho_t < -L) &\leq \sum_{x < y, y-x > L} \P\big(\;\eta^{x}_{t} \leq 0 < \eta^{y}_{t}, \; \hat T^{(x, y)} < t\;\big) \\
&\leq \sum_{z > L} \; \sum_{x \in \Z} \P\big(\;\eta^{x}_{t} \leq 0 < \eta^{x+z}_{t}, \; \hat T^{(x, x+z)} < t\;\big) \leq C \sum_{z \geq L} z e^{-cz},
\end{aligned}$$
which tends to zero as $L \to \infty$.
\eprthm

%%%%%%%%%%%%%% section 5 %%%%%%%%%%%%%%%%%%%%%%%%%%%%%%%%%%%%%%%%%%%%%%%%%%%%%%%%%%%%%%%%%%%%%%%%%%%%%%%%%%%%%%%%%%%%%%
\section{Estimate for a perturbed random walk}
In what follows, $\pi$ and $(\pi_z)_{z \in \Z}$ are probability distributions on $\Z$. We assume:
\be{symm} \pi \text{ is symmetric } (\textit{i.e. } \pi(-x) = \pi(x) \; \forall x);\ee
\be{support}\pi(x), \pi_z(x) > 0 \text{ for all } x \in \Z, z \in \Z - \{0\}; \ee
\be{pitail} \text{There exist } f, F > 0 \text{ such that } \pi(x), \pi_z(x) < F e^{-f|x|} \text{ for all } x \in \Z, z \in \Z; \ee
\be{totvar} \text{There exist } g, G > 0 \text{ such that } ||\pi_z - \pi||_{TV} < G e^{-g|z|} \text{ for all } z \in \Z.\ee
Given $x \in \Z$, let $\P_x$ be a probability under which a process $(X_n)$ is a Markov chain with transitions $P(z, w) = \pi_z(w-z)$ and $\P_x(X_0 = x)=1$. Define $H_0 = \inf\{n \geq 0: X_n = 0\}$.

\bt{rws} There exists $C > 0$ such that, for $x \in \Z$, $$\P_x(H_0 > N) < \frac{C|x|}{\sqrt{N}}.$$
\et

The proof of Theorem \ref{rws} will be carried out in a series of results. Fix $L >0$ such that $Ge^{-gL} < 1$ and let $I=[-L,L]$. Put $\epsilon_z = Ge^{-g|z|}$ for $z \in I^c$ and $\epsilon_z = 1$ for $z \in I$. A consequence of (\ref{totvar}) is that, for all $z \in \Z$, there exist probabilities $g_z, b_z^1, b_z^2$ on $\Z$ such that
\begin{eqnarray}
\pi_z = \epsilon_z b_z^1 + (1- \epsilon_z)g_z; \label{decomp1} \\
\pi = \epsilon_z b_z^2 + (1- \epsilon_z)g_z\label{decomp2}.
\end{eqnarray}
(Of course, if $z \in I$ we must have $b_z^1 = \pi_z, b_z^2 = \pi$).

We will construct the process $(X_n)$ coupled with other processes of interest. Let $(X_n, Z_n)$ be a Markov chain on $\Z \times \{0,1\}$ with transitions
\be{xztrans}Q((x, i), (y, j)) = \left\{
\begin{array}{ll} \epsilon_x \cdot b_x^1(y-x) &\text{ if } j = 1;\\ (1 - \epsilon_x) \cdot g_x(y-x) &\text{ if } j=0. \end{array} \right.\ee
We write $\P_x$ to represent any probability for this chain with $X_0 = x$, regardless of the law of $Z_0$. This abuse of notation is justified by the fact that $Z_0$ has no influence on the distribution of the other variables of the chain, and neither on the random variables to be defined below.
Let $(\mathcal{F}_n)$ be the natural filtration of the chain, and $T = \inf\{n \geq 1: Z_n = 1\}$.

Let $(\Psi_z)_{z \in \Z}$ be random variables defined on the same probability space as the chain above, independent of the chain and with laws $\Psi_z \stackrel{d}{=} b_z^2$. Additionally, let $(\Phi_n)_{n \geq 0}$ be a random walk with increment law $\pi$, initial state $0$, also defined on the same space as the previous variables and independent of them. For $n \geq 0$, define
\be{ycouple}
Y_n = \left\{ 
\begin{array}{ll}X_n, &\text{if } n < T;\\
X_{T-1} + \Psi_{X_{T-1}} + \Phi_{n - T}, &\text{if } n \geq T.\end{array} \right.
\ee
We can use (\ref{decomp1}) and (\ref{decomp2}) to check that under $\P_x$, $(X_n)$ is a Markov chain with transitions $P(z, w) = \pi_z(w-z)$ and initial state $x$, and $(Y_n)$ is a random walk with increment distribution $\pi$ and initial state $x$. They satisfy $X_n = Y_n$ on $\{T < n\}$.

Finally, we define some more stopping times. Let $H_0^Y = \inf\{n \geq 0: Y_n = 0\}, H_I = \inf\{n \geq 0: X_n \in I\}, \tau_0 = 0, \tau_1 = T \land H_I$ and $\tau_{k+1} = \tau_k + \tau_1 \circ \theta_{\tau_{k}}$ for $k \geq 1$, where $\theta_t$ denotes the shift operation $\theta_t((X_n, Z_n)_{n \geq 0}) = (X_{t+n}, Z_{t+n})_{n \geq 0}$. Note that $\tau_k \leq H_I$ for all $k$ and, if $\tau_k = H_I$ and $m > k$, then $\tau_m = H_I$. Also, $\tau_1 \leq H_0^Y$, because if $Y_n = 0$ for some $n$, then either $X_n = 0$, in which case $\tau_1 \leq H_I \leq H_0 \leq n$, or $X_n \neq 0$, in which case $\tau_1 \leq T < n$.

We will need the following standard facts about random walk on $\Z$:
\bl{reversible} 
$(i.)\; \P_x(H_0^Y > N) \leq \frac{C|x|}{\sqrt{N}}$ for some $C > 0$ and all $x \in \Z$;\\
$(ii.)\; \E_x(\#\{n < H_0^Y: Y_n = y\}) \leq C|y|$ for some $C > 0$ and all $x, y \in \Z$.
\el
\bpr $(i)$ is in \cite{spitzer}: see \textbf{P4} in Section 32 and Section 29. For $(ii)$, we have $\E_x(\#\{n < H_0^Y: Y_n = y\}) \leq \E_y(\#\{n < H_0^Y: Y_n = y\}) = \P_y(H^Y_0 < H^Y_{y+})^{-1}$, where $H^Y_{y+} = \inf\{n \geq 1: Y_n = y\}$, so it suffices to show that $\P_y(H^Y_0 < H^Y_{y+}) > c/y$ for some $c > 0$ and all $y \in \Z$. This can be done using Thomson's Principle for electric networks (see for example \cite{peres}, Theorem 9.10 and Section 21.2 for the infinite network case): if $y > 0$, take the unit flow $\theta(\overrightarrow{zw}) = 1$ if $z \in \{1, \ldots, y\}, w = z - 1$ and $\theta(\overrightarrow{zw}) = 0$ otherwise, and similarly if $y < 0$.
\epr

\bl{exptau} There exist constants $c, C > 0$ such that $\P_x\{|X_{\tau_1}| > r, \tau_1 < \infty, \tau_1 < H_I\} \leq Ce^{-cr}$ for all $x \in \Z, r \geq 0$.
\el
\bpr
If $x \in I$, then $\tau_1 = H_I = 0$ and the stated inequality is trivial. If $x \notin I$,
\begin{eqnarray*}
&&\P_x(|X_{\tau_1}| > r, \tau_1 < \infty, \tau_1 < H_I) =\\
&&\qquad\sum_{n=0}^\infty \sum_{z \in I^c} \P_x(X_0, \ldots, X_n \in I^c, X_n = z, Z_0 = \cdots = Z_n = 0, Z_{n+1} = 1, |X_{n+1}| > r) = \\
&&\qquad\sum_{n=0}^\infty \sum_{z \in I^c} \P_x(Y_0, \ldots, Y_n \in I^c, Y_n = z, Z_0 = \cdots = Z_n = 0) \cdot \P_z(Z_1 = 1, |X_1| > r) \leq\\
&&\qquad\sum_{n=0}^\infty \sum_{z \in I^c} \P_x(Y_0, \ldots, Y_n \in I^c, Y_n = z) \cdot \epsilon_z \cdot b_z^1\{w: |w| \geq |r-z|\} \leq\\
&&\qquad\sum_{z \in I^c} \inf\{\epsilon_z, \pi_z\{w: |w| \geq |r-z|\}\}\E_x(\#\{n < H_0^Y: Y_n = z\})\leq\\
&&\qquad C \sum_{z \in Z} |z|\inf\{Ge^{-g|z|}; Fe^{-f|r-z|}\} \leq Ce^{-cr}. 
\end{eqnarray*}
\epr 
\bc{corexptau} $(i.)\; A := \sup_{x \in \Z} \E_x(|X_{\tau_1}|; \tau_1 < \infty, \tau_1 < H_I) < \infty;$\\
$(ii.)\; \P_x(\tau_1 = \infty) = 0 \; \forall x\in \Z$.\\
Increasing $L$ if necessary,\\
$(iii.)\; \sigma := \inf_{x \in \Z} \P_x(\tau_1 = H_I < \infty)  > 0$;\\
$(iv.)\; \P_x(\tau_k < H_I) \leq (1-\sigma)^k$;\\
$(v.)\; \P_x(H_I = \infty) = 0$.
\ec
\bpr $(i.)$ is obtained by summing the two sides of the inequality of Lemma \ref{exptau} over $r$. For $(ii.)$, since $\tau_1 \leq H_0^Y$, $\P_x(\tau_1 > N) \leq \P_x(H_0^Y > N) \stackrel{N\rightarrow \infty}{\longrightarrow} 0$. For $(iii.)$, note that $$\begin{aligned}\P_x(\tau_1 = H_I < \infty) &= 1 - \P_x(\tau_1 = H_I = \infty) - \P_x(\tau_1 < \infty, \tau_1 < H_I)\\
&= 1 - 0 - \P_x(\tau_1 < \infty, |X_{\tau_1}| > L) \geq 1 - Ce^{-cL},\end{aligned}$$
which can be made positive by increasing $L$. Now, if $k \geq 1$,
$$\P_x(\tau_k < H_I) = \E_x(\mathds{1}_{\{\tau_{k-1} < H_I\}}\P_{X_{\tau_{k-1}}}(\tau_1 < H_I)) \leq (1-\sigma)\P_x(\tau_{k-1} < H_I)$$
by $(iii.)$, and continuing we get $(iv.)$ Finally, note that
$$\P_x(H_I=\infty) \leq \P_x(H_I=\infty, \tau_k < \infty \; \forall k) + \sum_{k=1}^\infty \P_x(H_I=\infty, \tau_k=\infty).$$
The first term is zero by $(iv.)$ and, using $(ii.)$,
$$\P_x(\tau_k = \infty) = \sum_{i=0}^{k-1} \P_x(\tau_i < \infty, \tau_{i+1} = \infty) = \sum_{i=0}^{k-1} \E_x(\mathds{1}_{\{\tau_i < \infty\}}\P_{X_{\tau_i}}(\tau_1 = \infty)) = 0,$$
so $(v.)$ follows.
\epr

\bl{controlT} There exists $C > 0$ such that, for all $x \in \Z$, $$\P_x(H_I > N) \leq \frac{C|x|}{\sqrt{N}}.$$
\el
\bpr
$$\begin{aligned}
\P_x(H_I > N) = &\P_x(H_I = \infty) + \sum_{k=0}^\infty\P_x(H_I > N, \tau_k < \tau_{k+1} = H_I < \infty)\\
&= \sum_{k=0}^\infty \P_x\left(\sum_{i=1}^{k+1}(\tau_i - \tau_{i-1}) > N, \tau_k < \tau_{k+1} = H_I < \infty\right)\\
&\leq \sum_{k=0}^\infty \sum_{i=1}^{k+1} \P_x\left(\tau_i - \tau_{i-1} > \frac{N}{k+1}, \tau_k < H_I\right)
\end{aligned}$$
We will show that, for $k \geq 0$ and $1 \leq i \leq k+1$,
\be{boundwithl1}
\P_x(\tau_i - \tau_{i-1} > l, \tau_k < H_I) \leq \frac{C|x|}{\sqrt{l}} (1-\sigma)^{k-2}
\ee
for some $C > 0$. So the above sum is less than
$$\sum_{k=0}^\infty \sum_{i=1}^{k+1} C |x| \sqrt{\frac{k+1}{t}}(1-\sigma)^{k-2} \leq \frac{C'|x|}{\sqrt{t}}$$
as required. To get (\ref{boundwithl1}), note that, if $i \leq k$, by Corollary \ref{corexptau}$(iv.)$, 
$$\P_x(\tau_i - \tau_{i-1} >l, \tau_k < H_I) = \E_x(\mathds{1}_{\{\tau_i - \tau_{i-1}> l\}} \P_{X_{\tau_i}}(\tau_{k-i} < H_I)) \leq (1- \sigma)^{k-i} \; \P_x(\tau_i - \tau_{i-1} > l),$$
so for any $i \in \{1, \ldots, k+1\}$,
\be{boundwithl2}
\P_x(\tau_i - \tau_{i-1} > l, \tau_k < H_I) \leq (1-\sigma)^{k-i} \;\P_x(\tau_i - \tau_{i-1} > l).
\ee
Now, using Lemma \ref{reversible},
$$\begin{aligned} \P_x(\tau_i - \tau_{i-1} > l) = \E_x(\mathds{1}_{\{\tau_{i-1} < H_I\}} \; \P_{X_{\tau_{i-1}}}(\tau_1 > l)) &\leq \E_x(\mathds{1}_{\{\tau_{i-1} < H_I\}} \; \P_{X_{\tau_{i-1}}}(H_0^Y > l))\\
&\leq (C/\sqrt{l})\E_x(\mathds{1}_{\{\tau_{i-1} < H_I\}} \;|X_{\tau_{i-1}}|).
\end{aligned}$$
If $i = 1$, the above expectation is equal to $|x|\mathds{1}_{\{x \in I^c\}}$; if $i>1$ it is equal to
$$\begin{aligned}
\E_x(\mathds{1}_{\{\tau_{i-2} < H_I\}} \; \E_x\left(|X_{\tau_{i-1}}|\cdot \mathds{1}_{\{\tau_{i-1} < H_I\}}|\mathcal{F}_{\tau_{i-2}})\right) &= \E_x\left(\mathds{1}_{\{\tau_{i-2} < H_I\}} \; \E_{X_{\tau_{i-2}}}(|X_{\tau_1}| \cdot \mathds{1}_{\{\tau_1 < H_I\}})\right)\\
&\leq A \P_x(\tau_{i-2} < H_I) \leq A (1-\sigma)^{i-2}
\end{aligned}$$
by Corollary \ref{corexptau} $(i.)$ and $(iv.)$. So, for any $i \in \{1, \ldots, k+1\}$,
\be{boundwithl3}
\P_x\{\tau_i - \tau_{i-1} > l\} \leq AC\frac{|x|}{\sqrt{l}}(1-\sigma)^{i-2}.
\ee
Putting together (\ref{boundwithl2}) and (\ref{boundwithl3}), we get (\ref{boundwithl1}).
\epr

From here to the proof of Theorem \ref{rws}, it is a matter of reapplying the ideas that established Corollary \ref{corexptau} and Lemma \ref{controlT}, so we simply sketch the main steps.

Define $T' = \inf\{n \geq 0: \{X_0, \ldots, X_n\} \cap I \neq \emptyset, \{X_0, \ldots, X_n\} \cap I^c \neq \emptyset\}, \lambda_0 = 0, \lambda_1 = T' \land H_0, \lambda_{k+1} = \lambda_k + \lambda_1 \circ \theta_{\lambda_k}$ for $k \geq 1$. From (\ref{support}), we get
\be{suppI1}\delta := \inf_{x \in I} P(x, 0) = \inf_{x \in I} \pi_x(-x) > 0. \ee
Two consequences are
\be{timeInI} \sup_{x \in I} \P_x(\lambda_1 > N) \leq (1-\delta)^N \ee
and
\be{suppI2} \inf_{x \in I} \P_x(\lambda_1 = H_0 < \infty) \geq \delta. \ee
Now, (\ref{timeInI}) and Lemma \ref{controlT} together imply
\be{lambdaInf} \forall x \in \Z, \P_x(\lambda_1 = \infty) = 0.\ee
Also, (\ref{suppI2}) gives 
\be{repetitions} \forall x \in \Z, \P_x(\lambda_k < H_0) \leq (1-\delta)^{\lfloor k/2 \rfloor}; \ee
this is justified by the fact that, if $\lambda_k < H_0$, then at least $\lfloor k/2 \rfloor$ times $X_n$ must have left $I$ without touching the origin. As in the proof of Corollary \ref{corexptau} $(v.)$, (\ref{lambdaInf}) and (\ref{repetitions}) are used to establish
\be{touch0}\forall x, \P_x(H_0 = \infty) = 0.\ee
The last ingredient is an analog of Corollary \ref{corexptau} $(i.)$,
\be{leavingI}B := \sup_{x \in I} \E_x(|X_{\lambda_1}|) = \sup_{x \in I} \E_x(|X_{\lambda_1}|; \lambda_1 < H_0) < \infty, \ee
which follows from (\ref{pitail}) and the fact that $I$ is finite.

We can now write
$$\P_x(H_0 > N) = \P_x(H_0 = \infty) + \sum_{k=0}^\infty \P_x(H_0 > N, \lambda_k < \lambda_{k+1} = H_0 < \infty)$$
and then, as in the preceeding proof, use (\ref{lambdaInf}), Lemma \ref{controlT}, (\ref{timeInI}), (\ref{repetitions}), and (\ref{leavingI}) to show that the above sum is less than $\frac{C|x|}{\sqrt{N}}$ for some $C > 0$.

To conclude, we mention the following result, for use in the proof of Lemma \ref{joinExp}. We omit its proof since it is simply a repetition of the above arguments.

\bl{discreteInv}
Let $H_{(-\infty, 0)} = \inf\{n: X_n < 0\}$. Then,
$$\sup_{x > 0} \E_x(\;|X_{H_{(-\infty, 0)}}|;\; H_{(-\infty, 0)} < H_0\;) \leq \sup_{x > 0} \E_x\; |X_{H_{(-\infty, 0)}}| < \infty.$$
\el

%%%%%%%%%%%% References %%%%%%%%%%%%%%%%%%%%%%%%%%%%%%%%%%%%%%%%%%%%%%%%%%%%%%%%%%%%%%%%%%%%%%%%%%%%%%%%%%%%%%%%%%%%%%%

\end{document}